\documentstyle[12pt]{amsart}

\input epsf
\epsfverbosetrue 

\newcommand{\xincludeps}[3]{\begin{figure}
                                \centerline{%
                                \epsfxsize=#2in
                                \epsfbox{#1}
                                }
                           \caption{ {\sf #3}  \label{fig:#1}}
                           \end{figure}
                          }

\newcommand{\yincludeps}[3]{\begin{figure}
                                \centerline{%
                                \epsfysize=#2in
                                \epsfbox{#1}
                                }
                           \caption{ {\sf #3}  \label{fig:#1}}
                           \end{figure}
                          }

\font\tenmsb=msbm10 scaled \magstep1
\font\sevenmsb=msbm8
\font\fivemsb=msbm6
\newfam\msbfam
\textfont\msbfam=\tenmsb
\scriptfont\msbfam=\sevenmsb
\scriptscriptfont\msbfam=\fivemsb
\def\Bbb{\fam\msbfam}

\newcommand{\rs}{\mbox{$\widehat{\Bbb C}$}}
\newcommand{\C}{\mbox{$\Bbb C$}}

\newcommand{\R}{\mbox{$\Bbb R$}}
\newcommand{\Z}{\mbox{$\Bbb Z$}}
\newcommand{\id}{\mbox{id}}
\newtheorem{thm}{Theorem}[section]
\newtheorem{defn}{Definition}

\newtheorem{prop}{Proposition}

\newtheorem{lemma}{Lemma}
\newtheorem{cor}{Corollary}

\newcommand{\vac}{\mbox{\O}}

\newcommand{\bdry}{\partial}
\newcommand{\cl}{\overline}
\newcommand{\gap}{\vspace{5pt}}
\newcommand{\rmk}{\gap \noindent {\bf Remark: }}
\newcommand{\intersect}{\cap}
\newcommand{\wt}{\widetilde}
\newcommand{\Int}{\mbox{Int}}
\newcommand{\Aut}{\mbox{Aut}}

\newcommand{\mult}{\mbox{mult}}
\newcommand{\be}{\begin{enumerate}}
\newcommand{\eb}{\end{enumerate}}
\newcommand{\bi}{\begin{itemize}}
\newcommand{\ib}{\end{itemize}}

\begin{document}

\author{Kevin M. Pilgrim}
\address{\hskip-\parindent Kevin M. Pilgrim\\
\\Mathematics Department\\Cornell University\\
       White Hall\\Ithaca, NY 14853\\USA}
\thanks{Research at MSRI is supported in part by NSF grant DMS-9022140}
\email{pilgrim@@math.cornell.edu}
\author{Tan Lei}
\address{\hskip-\parindent Tan Lei\\Department of Mathematics\\ 
       University of Warwick\\ Coventry CV4 7AL\\England}
\email{tanlei@@maths.warwick.ac.uk}

\title{Surgery on postcritically finite rational maps by blowing up an arc}

\begin{abstract}
Using Thurston's characterization of postcritically finite rational
functions as branched coverings of the sphere to itself, we give a new
method of constructing new conformal dynamical systems out of old
ones.  Let $f(z)$ be a rational map and suppose that the postcritical
set $P(f)$ is finite.  Let $\alpha$ be an embedded closed arc in the
sphere and suppose that $f|{\alpha}$ is a homeomorphism. Define a
branched covering $g$ as follows.  Cut the sphere open along $\alpha$.
Glue in a closed disc $D$.  Map $S^{2} - \Int (D)$ via $f$ and $\Int
(D)$ by a homeomorphism to the complement of $f(\alpha)$.  We prove
theorems which give combinatorial conditions on $f$ and $\alpha$ for
$g$ to be equivalent in the sense of Thurston to a rational map.  
The main idea in our proofs is a general theorem which forces 
a possible obstruction for $g$ away from the disc $D$ on which
the new dynamics is defined.
\end{abstract}

\maketitle

\section{Introduction}
\label{section:introduction}

A {\it rational map} $f(z)=p(z)/q(z)$ where $p$ and $q$ are relatively
prime complex polynomials determines a holomorphic map of the Riemann
sphere $\rs $ to itself, and so defines a holomorphic dynamical
system.  The set of points $C(f)$ where the derivative of $f$ vanishes
is the set of {\it critical points} of $f$; these are the points
where the local degree of the map is larger than one.  Counted with
multiplicity, there are $2d-2$ critical points where $d$ is the degree
of $f$.  The {\it postcritical set} $P(f)$ of $f$ is
defined by $P(f)=\cl{ \cup_{n>0} f^{n}(C(f)) }$.  The map $f$ is said
to be {\it postcritically finite} if $P(f)$ is finite.  The Riemann
sphere decomposes into the open {\it Fatou set} $F(f)$ of points whose
asymptotic behavior is stable under small perturbations, and the
closed {\it Julia set} $J(f)$ of points whose behavior under iteration
is chaotic.

A {\it branched covering} $f: S^{2} \to S^{2}$ is an
orientation-preserving continuous map such that for every $x \in
S^{2}$, there exist local coordinates near $x$ and $f(x)$ such that,
in these coordinates, $x$ and $f(x)$ are the point $0 \in \C$ and $f$
is given by $z \mapsto z^{n}$, $n \geq 1$.  The integer $n$ is called
the {\it local degree} of $f$ near $n$.  The sets $C(f)$ and $P(f)$ are
defined analogously as for rational maps.  

\gap\noindent {\bf Convention:}  We assume throughout this paper that
all rational maps and branched coverings are postcritically finite.

Two branched coverings $f$ and $g$ of the sphere to itself are said to
be {\it equivalent} if there are homeomorphisms $\psi_{0}$ and
$\psi_{1}$ such that $\psi_0 \circ f = g \circ \psi_1$, and such that
$\psi_0$ is isotopic to $\psi_1$ through homeomorphisms $\psi_t, t \in
[0,1] $ with $\psi_t | P(f) = \psi_0 | P(f)$ for all $t$.  

Thurston has given necessary and sufficient combinatorial conditions
for a postcritically finite branched covering $g$ to be equivalent to
a rational map, and has also shown that (with one understood
exception) this map is determined up to conformal congugacy by $g$;
see \cite{DH1}.  

If $f(z)$ is a rational map, suppose by cutting and pasting we
obtain a branched covering $g$.  When is $g$ equivalent to a rational
map? This question is, by Thurston's theorem, a combinatorial one.  In
general, given $d$ and $N$, there are infinitely many distinct
equivalence classes of branched coverings $f$ of degree $d$ with
$|P(f)|=N$.  Only finitely many of these can contain rational maps, so
it is useful to have constructions which, when applied to rational
maps, yield new rational maps.  In this paper we define a certain kind
of cutting and pasting process on $f(z)$ which is very robust: with
fairly weak hypotheses which do not include restrictions on the degree
of $f$ or its global combinatorics, the result $g$ will be equivalent
to a rational map.

The rough idea of our construction is the following.  Let $f(z)$ be a
postcritically finite rational map, let $\alpha$ be a closed, embedded
topological arc in $\rs $ and suppose that $f|\alpha$ is a
homeomorphism.  Define a new branched covering $g$ by cutting open the
sphere along $\alpha$, gluing in a disc $D$, mapping the complement of
$D$ by $f$, and mapping $D$ to the complement of $f(\alpha)$ in the
sphere.  A local model of $g$ near $D$ is the map $z \mapsto
\frac{1}{2}(z+1/z)$ on a neighborhood of the closed unit disc; this
maps the unit circle onto the interval $[-1,1]$.  Then
$\deg(g)=\deg(f)+1$, and the local degree of $g$ is one more than the
local degree of $f$ at the endpoints of $\alpha$.  More generally, one
can map the disc $D$ in a similar fashion $n$ times around the sphere
so that $\deg(g)=\deg(f)+n$.  We call this construction {\it blowing
up the arc $\alpha$ $n$ times}.  Note that by construction $f=g$ away
from $D$.  The precise conditions on $\alpha$ for this construction to
be well-defined, and the statements of our main theorems, are somewhat
technical and are given precisely in Sections \ref{subsection:conditions} 
and \ref{subsection:statements}.  We prove that under these conditions,
the combinatorial class of $g$ depends only on combinatorial data (Proposition  
\ref{prop:G well defined}).  

For example, we derive the following as a corollary to our first main
theorem, Theorem A.  To set up the statement, we recall that in a
periodic Fatou component of a postcritically finite rational map $f$,
there is a canonical family of arcs called {\it internal rays} which
are mapped homeomorphically under $f$ to other internal rays (see
Section \ref{subsection:construction of S} for the definition of
internal ray).

\gap
\noindent {\bf Theorem: }{\it If $\alpha$ is a periodic or preperiodic 
internal ray of $f$, then the branched covering $g$ which is $f$ blown
up $n$ times along $\alpha$ is equivalent to a rational map.}
\gap

We also apply the blowing up construction to obtain rational maps from
M\"{o}bius transformations of finite order (see Theorem 
\ref{subsection:blowing up moebius transformations}).  

The central idea in our proofs is Theorem 
\ref{thm:arcs intersecting obstructions}, 
{\it Arcs intersecting obstructions}.  This theorem asserts that under
suitable hypotheses, an obstruction $\Gamma$ to the existence of a
rational map equivalent to a branched covering $g$ cannot intersect
any family of arcs $\Lambda$ which maps in a certain way.  We use this
theorem to ``push'' an obstruction away from the region $D$ in the
sphere where we add new dynamics.  We then conclude that the
obstruction $\Gamma$ yields an obstruction for the original rational
map $f$, which is impossible.

Our theorems may also be thought of as {\it combination theorems}
between conformal dynamical systems.  Other examples of such theorems
include {\it mating} (see \cite{tanlei:matings},
\cite{shishikura:tan:cubicmating}, \cite{jluo:thesis} ), where two
polynomials of the same degree are glued together along their circles
at infinity, and {\it tuning} (see \cite{rees:parami},
\cite{ahmadi:thesis}), where the dynamics of a polynomial is glued
into the dynamics of a rational map near a superattracting cycle.  In
both of these examples the degree of the new map is the same as that
of the old one.  To date, all known combination theorems involving
mating and tuning have hypotheses on the degree of the map and on the
combinatorics of both the original map $f(z)$ and the new dynamics
which is being glued into $f(z)$.  Also, a definition of mating or
tuning for branched coverings has not yet been given in such a manner
as to be defined in terms of purely combinatorial data.  Our
techniques of proof also allow us to establish that certain blown-up
quadratic polynomials may be mated (in a suitably generalized way)
with an {\it arbitrary} quadratic polynomial to yield a rational map;
see Section \ref{subsection:generalized matings}.  Moreover, the
technique generalizes easily to similar statements in higher degrees.

\gap \noindent {\bf Organization of the paper.}  

\begin{itemize}

\item \S \ref{section:theorems}: we state needed definitions, give
the precise definition of the blowing up construction, state our first
two main theorems (Theorems A and B) precisely, and prove that the
combinatorial class of $g$ depends only on combinatorial data
(Proposition \ref{prop:G well defined}).  

\item \S \ref{section:proofs}: we prove Theorems A and B.

\item \S \ref{section:complements}: we prove other related results.

\item \S \ref{section:examples}: we give examples illustrating 
our main theorems.  

\item \S \ref{appendix:arcs}: is an appendix containing a proof of 
Theorem \ref{thm:arcs intersecting obstructions}, {\it Arcs intersecting
obstructions}.  

\end{itemize}

\noindent {\bf Acknowledgments:}  The pictures were drawn using programs
written by Curt McMullen.  The authors also thank MSRI, where some 
of this work was prepared.

\gap\noindent {\bf Example: blowing up a preperiodic arc in the
basilica.} Let $f(z)=z^{2}-1$.  Let $\alpha = [0,1]$. The filled-in
Julia set of $f(z)$ is the set of points whose orbits are bounded;
this set is the black region in Figure \ref{fig:basilica.ps}.  Let $g$
be the branched covering which is $f$ blown up once along $\alpha$.
The Julia set of the rational map equivalent to $g$ is shown in Figure
\ref{fig:semibasilica.ps}; see also \cite{kmp:jcurve} for a discussion of
this example and for a description of $g$ as a branched covering. 

\yincludeps{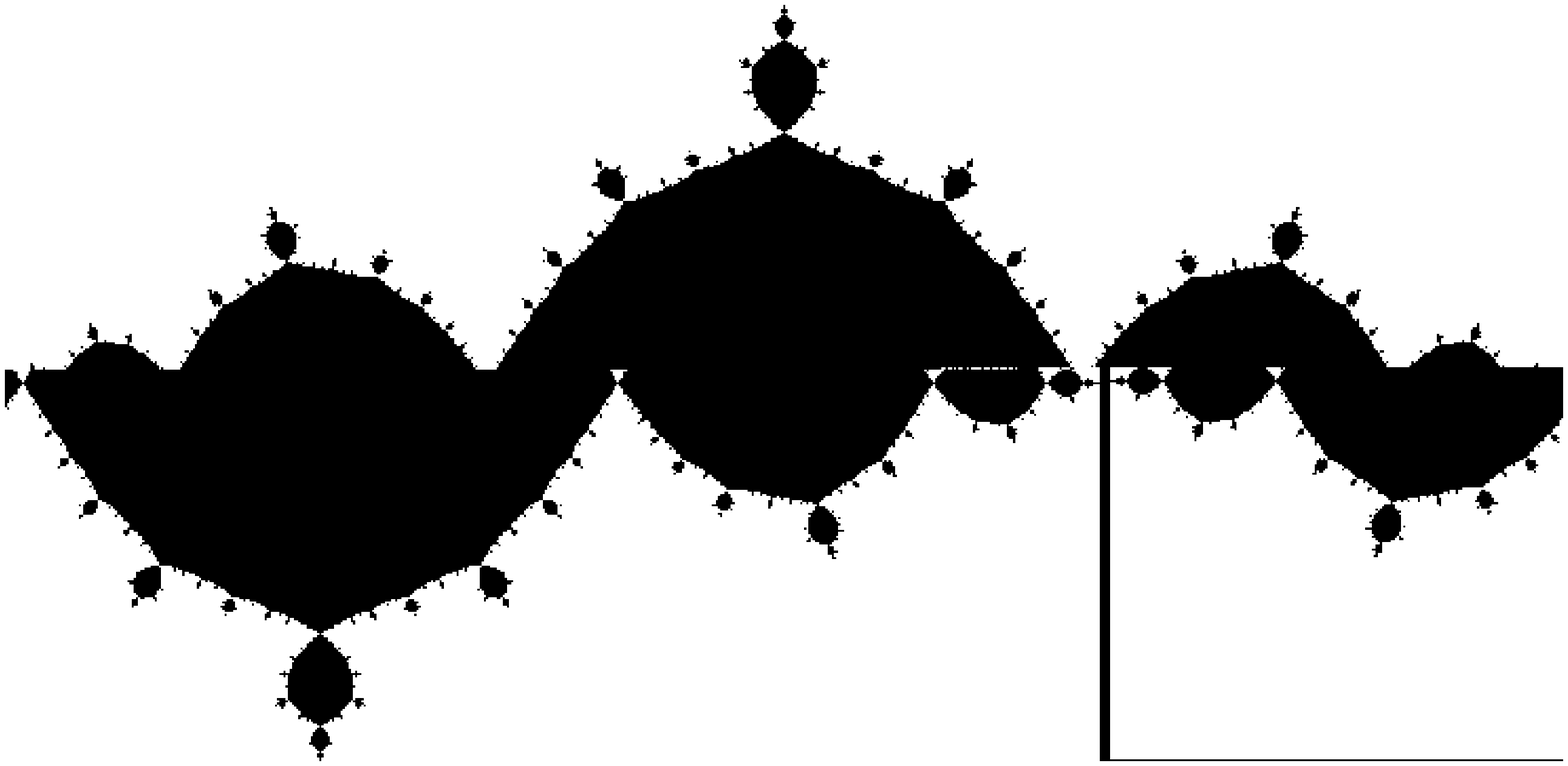}{1.75}{The filled-in Julia set of $f(z)=z^2-1$.}

\yincludeps{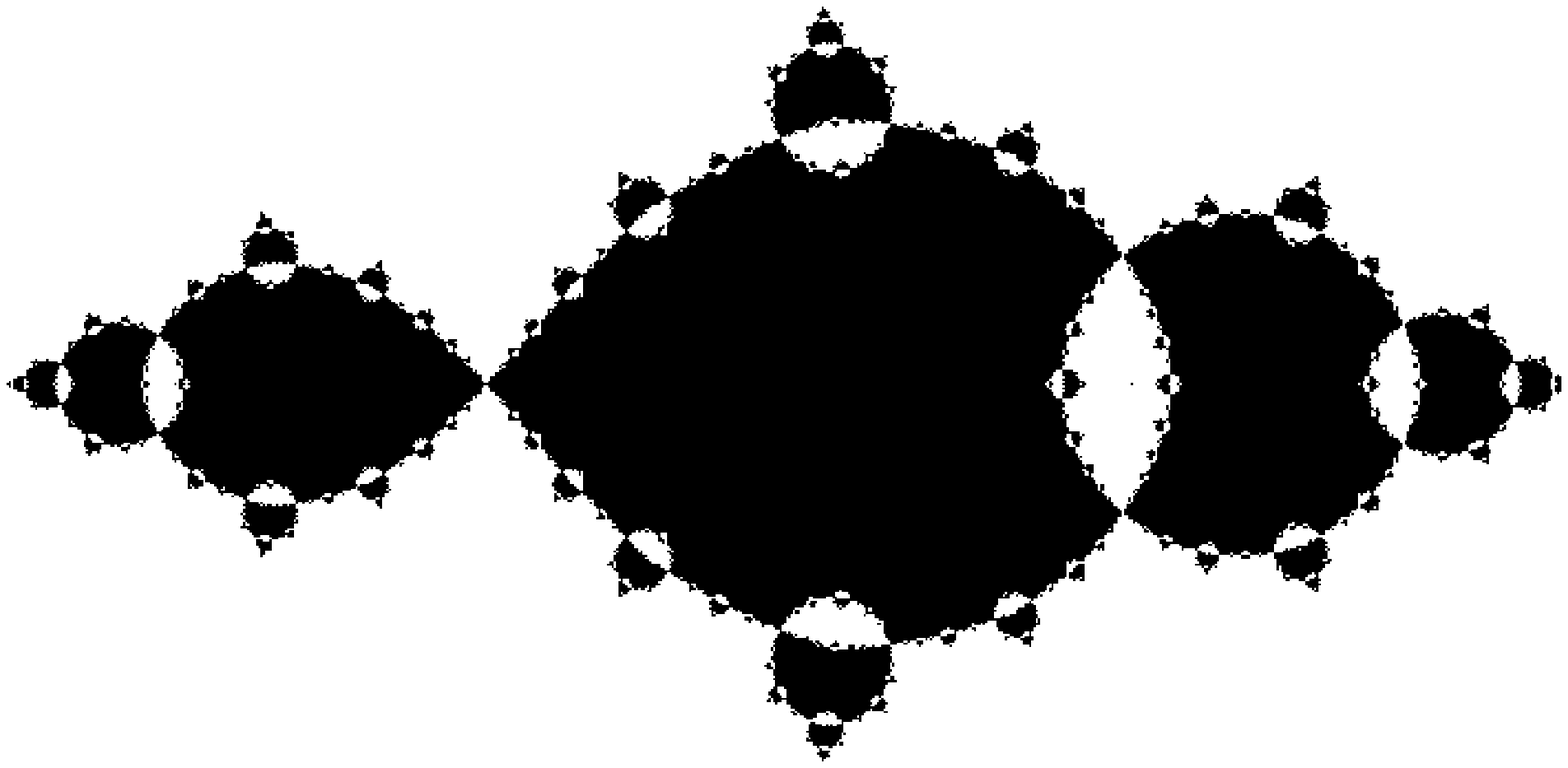}{1.75}{Julia set of the rational map corresponding 
to blowing up $f(z)=z^{2}-1$ along $[0,1]$.  The white regions form the basin of
infinity.}

\section{Blowing up an arc}
\label{section:theorems}

\S\S 2.1 and 2.2 give definitions needed to formulate our construction
and theorems precisely.  \S 2.3 states conditions under which blowing
up will be a well-defined operation on equivalence classes.  \S 2.4
gives the precise statements of Theorems A and B.  \S 2.5 gives the
construction of blowing up an arc, a proof that the equivalence class
of the new map $g$ depends only on combinatorial information, and a
proposition ({\it Orbits off $\alpha$ persist}) which we will use in
our proofs.

\subsection{Marked branched covers} 

\begin{defn}  A {\bf marked branched covering} is a pair $(f,X)$ where
$f: S^2 \to S^2$ is a branched covering and $X$ is a finite set
containing $P(f)$ such that $f(X) \subset X$.  
\end{defn}

Two marked branched coverings $(f,X), (g,Y)$ are said to be {\it
equivalent} if they satisfy the definition of equivalence of branched
coverings with the sets $P(f)$ and $P(g)$ replaced by $X$ and $Y$
respectively.  There is a functor from the category of marked branched
coverings $(f,X)$ and equivalences to the category of branched
coverings and equivalences. The image of a marked brached covering
$(f,X)$ under this functor is $f$; the image of an equivalence
$\psi_{0}:(S^{2},X) \to (S^{2},Y)$ is the map of pairs
$\psi_{0}:(S^{2},P(f)) \to (S^{2},P(g))$.  This functor ``forgets''
that $\psi_{0}$ must send $X-P(f)$ to $Y-P(g)$. We call this the {\it
forgetful functor} sending $(f,X)$ to $f$.

\gap \noindent {\bf Examples:}  
\begin{enumerate}
\item Let $f$ be a rational map and $X$ a finite forward-invariant set
containing $P(f)$.  Then $(f,X)$ is a branched cover marked by the set
$X$.  We will call $(f,X)$ a {\it marked} rational map.  Note that if
$(f,X)$ is a marked rational map with $\deg(f) = 1$ and $|X| \geq 3$,
then $f$ is necessarily a M\"{o}bius transformation of finite order.  

\item Let $p(z)=z^{2}-1$.  Then $p:[-1,0] \to [-1,0]$ is an
orientation-reversing homeomorphism.  Postcompose $p(z)$ with a
homeomorphism $h$ isotopic to the identity relative to $P(p)$ such
that for $f(z)=h \circ p(z)$, $f([-1/4,0])=[-1,-3/4]$ and $f([-1,-3/4]) =
[-1/4,0]$.  Let $X=P(f) \cup \{-1/4, -3/4\}$.  Then $(f,X)$ is a
branched cover marked by the set $X$.

\end{enumerate}

In both examples, the image of $(f,X)$ under the forgetful functor is
combinatorially equivalent to a rational map.  In the first case,
however, the set $X-P(f)$ of ``new'' marked points are already
preperiodic under the dynamics of a rational map, whereas in
the second case we have added new preperiodic points by our choice of
representative of the combinatorial class of $p(z)$.

\subsection{Arcs and curves in $(S^{2},X)$}

Let $X$ be a finite set in the sphere.  An {\it arc $\alpha$ in
$(S^{2},X)$} we define to be the image of the unit interval $[0,1]$ under
a map $j$ into the sphere such that the endpoints $e(\alpha):=j(\{0,1\})
\subset X$, $j|(0,1)$ is an embedding, and $j(0,1) \cap X= \vac$. Thus
we allow the endpoints $e(\alpha)$ of $\alpha$ to coincide (though we
will not cut along such arcs in our construction).  A {\it simple
closed curve in $(S^{2},X)$} is a simple closed curve in $S^{2}-X$; it
is said to be {\it peripheral} if it is isotopic into every
neighborhood of some point in $X$ and {\it essential} if it is 
not contractible in $S^2-X$.  

Two arcs (resp. curves) $\eta_{0}, \eta_{1}$ in $(S^2,X)$ are said to
be {\it isotopic relative to $X$}, written $\eta_{0} \simeq_X
\eta_{1}$, if there is a continuous, one-parameter family $\eta_{t}, t
\in [0,1]$ of such arcs (respectively curves) joining $\eta_{0}$ to
$\eta_{1}$.  The isotopy class relative to $X$ of an arc or simple
closed curve $\eta$ we will denote by $[\eta]_X$, or by $[\eta]$ when
the set $X$ is understood.

A set of pairwise nonisotopic arcs (essential curves) in $(S^2,X)$
will be called an {\it arc (curve) system} in $(S^2,X)$.  If $\Lambda$
and $\Lambda'$ are arc or curve systems in $(S^2,X)$ we write $\Lambda
\simeq_X \Lambda'$ if every element of $\Lambda$ is isotopic to a
unique element of $\Lambda'$ relative to $X$, and conversely.

The {\it intersection number} between a pair $\alpha, \beta$
consisting of either arcs or curves in $(S^2,X)$ is
defined by 

\[ \alpha \cdot \beta = \min_{\alpha \simeq_X \alpha', \beta \simeq_X \beta'} 
\# \{(\alpha' - e(\alpha')) \cap (\beta' - e(\beta'))\}.\]

\gap \noindent {\bf Invariance up to isotopy of arc systems.}

Let $(f,X)$ be a marked branched covering.

If $\lambda$ is an arc in $(S^2,X)$, we will call the closure
$\wt{\lambda}$ of a connected component of $f^{-1}(\lambda -
e(\lambda))$ a {\it lift} of $\lambda$ under $f$.  Thus each arc
$\lambda$ in $(S^2,X)$ has $\deg(f)$ distinct lifts $\wt{\lambda}$.
If the endpoints of $\lambda$ are distinct, each lift maps
homeomorphically under $f$ onto $\lambda$. If $\Lambda$ is an arc
system in $(S^2,X)$, a {\it lift} $\wt{\Lambda}$ is defined to be an
arc system in $(S^2,X)$ each of whose elements $\wt{\lambda}$ is a
lift of some element $\lambda \in \Lambda$.

An arc system $\Lambda$ in $(S^2,X)$ is said to be {\it forward
invariant under $f$ up to isotopy relative to $X$} if there exists a
subset $\Lambda_0 \subset \Lambda$ (possibly all of $\Lambda$) and a
lift $\wt{\Lambda}_0$ of $\Lambda_0$ such that $\wt{\Lambda}_0
\simeq_X \Lambda$. If $\Lambda \simeq_X \Lambda'$ and $\Lambda$ is
forward invariant up to isotopy relative to $X$, then so is
$\Lambda'$.  (An isotopy between an element $\lambda \in \Lambda$ and
an element $\lambda' \in \Lambda'$ can be lifted to an isotopy between
$\wt{\lambda} \in \wt{\Lambda}$ and an arc $\wt{\lambda}'$ which is a
lift of $\wt{\lambda}'$.)

The next proposition allows us to extract arc systems which are
forward-invariant up to isotopy from a collection of arcs which are
literally forward-invariant.

\begin{prop}
\label{prop:finding invariant arc systems}
Suppose $\cal{L}$ is a collection of arcs $l$ in $(S^2,X)$ such that
for each $l \in \cal{L}$, $f|l$ is a homeomorphism and $f(l) \in
\cal{L}$.  Then there is a subset $\Lambda \subset \cal{L}$ 
which is an arc system forward invariant under $f$ up to isotopy, and
such that each $l \in \cal{L}$ is isotopic relative to $X$ to a unique
element of $\Lambda$.
\end{prop}

\pf Let $\Lambda_0 \subset f(\cal{L})$ be an arc system chosen
arbitrarily so that each arc $f(l), l \in \cal{L}$ is isotopic
relative to $X$ to a unique element of $\Lambda_0$.  Then for each
$\lambda \in \cal{L}$, $f(\lambda) \simeq_X
\lambda_0$ for a unique $\lambda_0 \in \Lambda_0$.  Hence by lifting
isotopies, $\lambda \simeq_X \wt{\lambda}_0$, where $\wt{\lambda}_0$
is a lift of $\lambda_0$.  Hence every element of $\cal{L}$ is
isotopic to a lift of an element of $\Lambda_0$.  We now enlarge our
set $\Lambda_0$ so that the last property holds.  Let $\cal{L}$$_1$
denote the set of elements of $\cal{L}$ which are not isotopic to any
element of $\Lambda_0$.  Let $\Lambda_1 \subset $$\cal{L}$$_1$ be any
arc system for which each element of $\cal{L}$$_1$ is isotopic to a
unique element of $\Lambda_1$, and for which each element $\lambda_1
\in \Lambda_1$ is isotopic to some element of $\cal{L}$$_1$.  Then
$\Lambda = \Lambda_0 \cup \Lambda_1$ is an arc system with the desired
properties.
\qed

\subsection{Conditions under which blowing up is defined}
\label{subsection:conditions}

Now suppose that $(f,X)$ is a branched covering of degree $\geq 1$
marked by the finite set $X$, and let $\alpha$ be an arc in $S^{2}$.
Assume that the following conditions are satisfied:
\gap

\noindent {\bf Blowing Up Conditions:}
\begin{enumerate}
\item $f|\alpha$ is a homeomorphism; 
\item $\alpha$ is a union of arcs $\alpha_{j},
j=1,2,...,L$ in $(S^{2},f^{-1}(X))$;
\item $\Int(\alpha) \cap (X  \cup C(f))= \vac$.
\end{enumerate}

Note that the arcs $f(\alpha_j)$ are arcs in $(S^2,X)$.  

\gap\noindent {\bf Remarks:}  
Allowing $\alpha$ to consist of more than one arc in
$(S^{2},f^{-1}(X))$ is indeed useful: for example, consider the
example of blowing up the ``airplane'' in Section
\ref{section:examples}.   More generally, we may apply the blowing 
up construction to finite collections $\alpha^1, \alpha^2, ...$ of
arcs for which $\Int(\alpha^i) \cap \alpha^j = \vac$ when $i \neq j$
and for which each arc satisfies the Blowing Up Conditions.

\subsection{Statement of main results}
\label{subsection:statements}

We assume the following in this subsection. Let $(f,X)$ be a marked
rational map with $\deg(f) \geq 2$.  Let $\alpha$ be an arc in $\rs$
satisfying the Blowing Up Conditions.  Let $(g,X)$ be marked branched
cover which is the map $(f,X)$ blown up $n$ times along $\alpha$.

\gap\noindent{\bf Theorem A}
{\it If for all $j=1,2,...,L$, 
\be

\item $f(\alpha_j)$ is isotopic relative to $X$ to an arc which is
contained in a finite arc system $\Lambda_j$ which is forward
invariant under $f$, and
 
\item there is a subset $\Lambda_{j,0} \subset \Lambda_j$ with a lift
$\wt{\Lambda}_{j,0}$ such that $\wt{\Lambda}_{j,0} \simeq_X \Lambda_j$
and $\wt{\Lambda}_{j,0} \cap \Int(\alpha) = \vac$,  

\eb
then $g$ is combinatorially equivalent to a rational map.}

\begin{cor}
If $\alpha$ is a periodic or preperiodic internal ray of $f$, then $g$
is equivalent to a rational map.
\end{cor}

\begin{cor}
Suppose $\alpha$ is a finite union of periodic or preperiodic internal rays
of $f$.  If for all  $n \geq 0$, the intersection number relative
to $X$ satisfies $f^n (f(\alpha)) \cdot f(\alpha) \neq 0$ only when 
$f^n (f(\alpha)) = f(\alpha)$, then $g$ is equivalent to a rational map. 
\end{cor}

The definition of internal ray is given in Section
\ref{subsection:construction of S}. 

\gap Next, assume $f(z)$ is a polynomial.  According to work of Douady 
and Hubbard, there is a canonically defined finite topological tree
$T$ contained in the filled-in Julia set of $f$, depending only on $f$
and $e(\alpha)$, such that $f(T) \subset T$ with vertices mapped to
vertices. In Section \ref{subsection:proof of theorem two} we outline
the construction of $T$.

\gap\noindent{\bf Theorem B: (Blowing up an arc in a generalized
Hubbard tree)} 
{\it If $\alpha \subset T$, then $g$ is combinatorially equivalent to
a rational map.}

\rmk The Blowing Up Conditions strongly restrict the possibilities for
$\alpha$.  The novel feature of this theorem, relative to the first,
is that we do not require that the $\alpha_j$ are eventually periodic 
up to isotopy relative to $X$.

\subsection{Definition of blown-up map}

Let $(f,X)$ and $\alpha$ satisfy the Blowing Up Conditions of Section
2.3 .  The construction is outlined in Figure \ref{fig:blowup.ps} for the
case $n=4$.  

\xincludeps{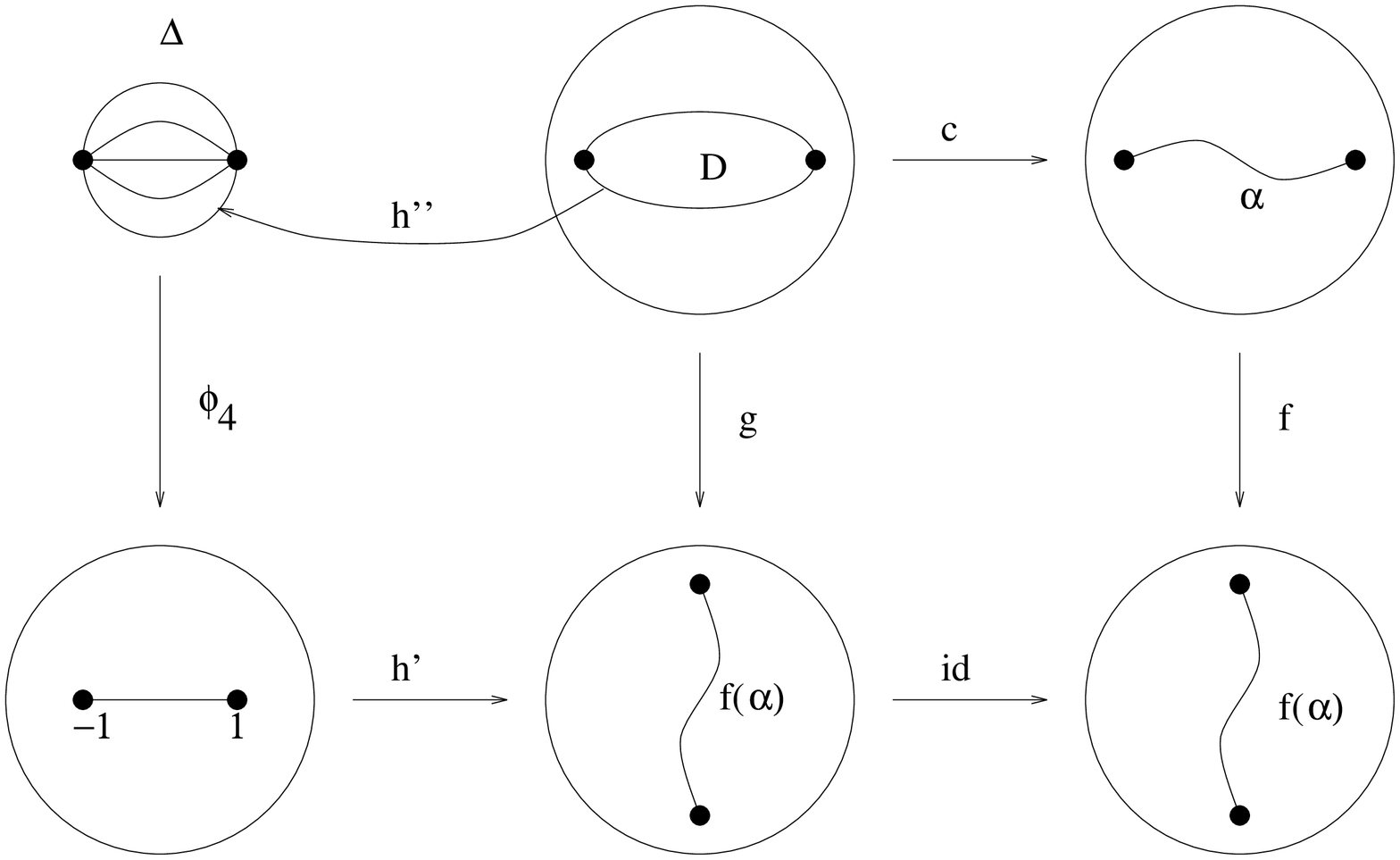}{4.5}{Construction of the blown-up map $g$ from $f$.}
We first cut the sphere along the interior of $\alpha$.  Let $W$ be an
open disc containing $\Int(\alpha)$ bounded by a Jordan curve such that
$\cl{W} \intersect f^{-1}(X) = \alpha \intersect f^{-1}(X)$ and $\bdry
W \supset e(\alpha)$.  Since $\Int(\alpha) \cap(X \cup C(f)) = \vac$,
we may choose $W$ so that $f|\cl{W}$ is a homeomorphism.  Let $h_{t}: \rs -
\Int(\alpha) \hookrightarrow \rs $ be a continuous family of
embeddings satisfying 

\begin{enumerate}

\item $h_{0}=\id $;

\item $h_{t}|(\rs - W) = \id $ for all $t$;

\item $h_{1}(\rs - \Int(\alpha)) = \rs - (D-e(\alpha))$, where $D$ is
a closed disc.

\end{enumerate}
The Jordan curve $\bdry D$ is a union of two closed arcs
$\bdry D_{+}$ and $\bdry D_{-}$ with endpoints $e(\alpha)$ which we
denote by $\{e_{-1},e_{+1}\}$.  The map $h_{1}^{-1}$ extends uniquely
to a continuous ``collapsing'' map $c: \rs - \Int(D) \to \rs $ sending
$\rs - (D-e(\alpha))$ homeomorphically to $\rs - \Int(\alpha)$ and
each arc $\bdry D_{+}$ and $\bdry D_{-}$ homeomorphically to $\alpha$.

We add the following new dynamics.  Let $\Delta = \{z\; |\; |z|<1\; \}$,
and let $\phi_{n}: \cl{\Delta} \to \rs $ be the map given by
\[ \phi_{n}(z)=A_{2}( A_{1}(z) )^{2n}, 	\]
where $A_{1}(z)=-i \frac{z+1}{z-1}$ and $A_{2}(z)=\frac{z-1}{z+1}$.
The M\"{o}bius transformation $A_{1}$ is a rigid rotation with respect
to the spherical metric sending the ordered triple $(-1,0,1)$ to
$(0,i,\infty)$ and sending $\cl{\Delta}$ to $Im(z) \geq 0$; the map $z
\mapsto z^{2n}$ maps the half-plane $Im(z) \geq 0$ $n$ times around
the sphere sending the real axis to the nonnegative real axis; the
M\"{o}bius transformation $A_{2}$ is a rigid rotation sending
$(0,1,\infty)$ to $(-1,0,1)$.  Thus the map $\phi_{n}$ wraps the
closed unit disc $\cl{\Delta}$ around the sphere $n$ times in such a way
that $\phi_{n}$ is locally injective on $\cl{\Delta}-\{-1,1\}$,
$\phi_{n}(-1)=-1, \phi_{n}(1)=1$, and $\phi_{n}(\bdry \Delta)=[-1,1]$.
Let $\bdry \Delta_{+}, \bdry
\Delta_{-}$ denote respectively the intersection of $\cl{\Delta}$ with
$Im(z) \geq 0$ and $Im(z) \leq 0$.

We now describe the gluing of $f$ and $\phi_{n}$.  Choose a
homeomorphism $$h': (\rs , [-1,1]) \to (\rs , f(\alpha))$$ which sends
$-1$ to $f(e_{-1})$ and $+1$ to $f(e_{+1})$.  The compositions
$(\phi_{n}|_{\bdry \Delta_{\pm}})^{-1} \circ (h')^{-1} \circ f \circ
c |_{\bdry D_{\pm}}: \bdry D_{\pm} \rightarrow \bdry \Delta_{\pm}$
are homeomorphisms of closed arcs and together define a homeomorphism
from $\bdry D$ to $\bdry \Delta$.  Since $D$ and $\cl{\Delta}$ are
closed discs we can extend this to a homeomorphism $h''$ between $D$
and $\cl{\Delta}$.

We now define a new branched covering $g$ by 
\[ g(z) =
\left\{ 
\begin{array}{ll}
   f  \circ c (z),		  &  \mbox{if $ z \in \rs - \Int(D)$}, \\
   h' \circ \phi_{n}\circ h''(z), &  \mbox{if $ z \in \cl{D}.$} 
\end{array}
\right. 
\]
The two definitions agree on the boundary, by construction, hence they
give a well-defined branched covering of the sphere to itself such
that $\deg(g)=\deg(f)+n$.   

We may also apply the above construction simulaneously to finite
collections of arcs $\alpha^1, \alpha^2, ...$ which satisfy the
Blowing Up Conditions and which satisfy $\Int(\alpha^i) \cap
\alpha^j = \vac, i \neq j$.  We omit the details of this more general
construction but emphasize that we may assume that for $i \neq j$,
$\cl{W^i} \cap \cl{W^j} \subset e(\alpha^i) \cap e(\alpha^j) \subset X$.  

\begin{prop}[$g$ depends only on classes]
\label{prop:G well defined}
The combinatorial class of $g$ depends only on the integer $n$, 
the combinatorial class of $(f,X)$, and the isotopy classes of
the arcs in $(S^{2},f^{-1}(X))$ comprising $\alpha$.  
\end{prop}

\pf We prove in detail only the case of blowing up along a single arc;
the case of blowing up along multiple arcs in similar, since we may
choose the closed discs $\cl{W^i}$ to intersect only in $X$.  

If $(f_{1},X_{1})$ is equivalent to $(f_{2}, X_{2})$, then
$(f_{1},f_{1}^{-1}(X_{1}))$ is equivalent to $(f_{2},
f_{2}^{-1}(X_{2}))$: an equivalence $\psi_{0}: (S^{2},X_{1}) \to
(S^{2},X_{2})$ may be lifted to an equivalence
$\psi_{1}:(S^{2},f_{1}^{-1}(X_{1})) \to (S^{2}, f_{2}^{-1}(X_{2}))$.

Let $\phi_{0}, \phi_{1}: (S^{2},f_{1}^{-1}(X_{1})) \to
(S^{2},f_{2}^{-1}(X_{2}))$ be a pair of homeomorphisms giving an
equivalence between $(f_{1},f_{1}^{-1}(X_{1}))$ and
$(f_{2},f_{2}^{-1}(X_{2}))$.  Then $\phi_{0} \circ f_{1}=f_{2}\circ
\phi_{1}$ and $\phi_{0} \simeq \phi_{1}$ relative to
$f_{1}^{-1}(X_{1})$.  Let $\alpha^{2}=\phi_{1}(\alpha^{1})$.  Then
$\alpha^{2}$ is an arc in $S^{2}$ which is comprised of arcs in
$(S^{2}, f_{2}^{-1}(X_{2}))$.  We must prove that if $(g_{i},X_i)$ is
$(f_{i},X_{i})$ blown up $n$ times along $\alpha^{i}, i=1,2$, then
$g_{1}$ and $g_{2}$ are combinatorially equivalent.  For convenience,
if $\beta^i, i=1,2$ are arcs which are unions of finitely many arcs
$\beta^i _j$ in $(S^2,Y)$, we will write $\beta^1 \simeq_Y \beta^2$ if
$\beta^1 _j \simeq_Y \beta^2 _j$ for all $j$.

\gap\noindent{\bf Step 1:}  We first show the following.  Suppose for
$i = 1,2$, $g^i$ is the map produced from the map $f$ by blowing up an
arc $\alpha^i$, where $\alpha^1 \simeq_{f^{-1}X} \alpha^2$, and $W_i,
D_i, h'_i, h''_i$, and $h^i _t$ are the discs and maps used in the
definition of $g_i$.  Then $g_1$ is equivalent to $g_2$.

Since $f|\cl{W}_i$ is a homeomorphism, $f(\cl{W}_i)$ is a closed disc.
There exists a homeomorphism $\psi_0: (S^2,X) \to (S^2,X)$ isotopic to
the identity relative to $X$ such that $\psi_0 (f(\cl{W}_1),
f(\alpha^1)) = (f(\cl{W}_2), f(\alpha^2))$.  For since $\alpha^1
\simeq_{f^{-1}X}\alpha^2$, $f(\alpha^1) \simeq_X f(\alpha^2)$.  By
results of Epstein \cite{dbaepstein:curves}, there is an ambient
isotopy of $S^2$ carrying $\bdry f(\cl{W}_1) \cup f(\alpha^1)$ to
$\bdry(f(\cl{W}_2)) \cup f(\alpha^2)$ fixing $X$ pointwise.  Let
$\Psi: S^2 \times I \to S^2$ be such an isotopy such that $\Psi (\cdot, 0)
= \id$, $\Psi(\cdot,1) = \psi_0$, and $\Psi(X,t) = \id_X$ for all $t$.

Since $X \supset P(f)$, $f| S^2 - f^{-1}X \to S^2 - X$ is a covering
map.  Hence there is a lift $\wt{\Psi}$ of $\Psi$ which is an
isotopy fixing $X$ pointwise between the identity and a map which we
call $\psi^1$.  Note that $\psi^1: (\cl{W}_1, \alpha^1) \to (\cl{W}_2,
\alpha^2)$.  

Recall that $g_i | S^2 - \cl{W}_i = f| S^2 - \cl{W}_i, i=1,2$.  We
will define a homeomorphism $\psi_1: (S^2,X) \to (S^2,X)$ such that
$g_2 \psi_1 = \psi_0 g_1$ and show $\psi_0$ is isotopic to $\psi_1$
relative to $X$.  The map $\psi_1$ will be defined as four maps
$\psi^{(1)}=\psi^1 | S^2 - W_1$, $ \psi^{(2)},
\psi^{(3)}, \psi^{(4)}$ which piece together.

We now define $\psi^{(2)}$.  Let $K_i = (g_i | D_i)^{-1}(f(\alpha^i)),
i=1,2$.  Then $K_i$ is a graph in $S^2$ with exactly two vertices
which are the endpoints $e(\alpha^i) \subset f^{-1}(X)$ and with $n+2$
edges, each of which joins one vertex to the other.  There is a
well-ordering on the set of edges of $K_i$: $E < E'$ if the edge $E$
separates $\Int(E')$ from $\Int(\bdry_{-}D_i)$ in $D_i$.  Write the
set of edges in increasing order as $\{E^i _0, ..., E^i _{n}\}$.
Since we glue in the same new dynamics $\phi_n$ in the construction of
$g_1$ and $g_2$, there is a homeomorphism $\psi^{(2)}: K_1 \to K_2$
such that $g_2 \psi^{(2)} = \psi_0 g_1 $ and $\psi^{(2)}(E^1 _l) = E^2
_l, l=0,1,...,n$, i.e. $\psi^{(2)}$ preserves the well-ordering of the
edges.

We now define $\psi^{(3)}$.  For $l=0,1,...,n-1$, the complementary
component $U^i _l$ of $D^i - K_i$ bounded by $E^i _l \cup E^i _{l+1}$
is a Jordan domain mapping homeomorphically under $g_i$ onto $S^2 -
f(\alpha^i)$, $i=1,2$.  Hence for each $l$, there is a homeomorphism
$\psi^{(3),l}: U^1 _l \to U^2 _l$ such that $g_2 \psi^{(3),l} | U^1 _l
= \psi_0 g_1 | U^1 _l$.  Since $\psi^{(2)}$ preserves the ordering of
edges, the homeomorphisms $\psi^{(3),l}$ and $\psi^{(2)}$ agree on the
common boundaries of their domains of definition.

A similar argument shows that there is a homeomorphism $\psi^{(4)}:
W_1 - D_1 \to W_2 - D_2$ such that $g_2 \psi^{(4)}| W_1 - D_1 = \psi_0
g_1 | W_1 - D_1$ and which agrees with $\psi^{(1)}$ on $\bdry W_1$ and
with $\psi^{(2)}$ on $\bdry D_1$.  

Hence the maps $\psi^{(1)}, \psi^{(2)}, \psi^{(3)}, \psi^{(4)}$ define
a homeomorphism of the sphere $\psi_1$ such that $\psi_0 g_1 = g_2
\psi_1$.  

That $\psi_0$ and $\psi_1$ are isotopic relative to $X$ follows from
the fact that $\psi_1 | S^2 - W_1$ is isotopic (via $\wt{\Psi}$) to
the identity relative to $X$ and the fact that $X \cap \cl{W}_1
\subset e(\alpha^1) \subset \bdry W_1$.  This completes the proof of
 Step 1.

\gap\noindent{\bf Step 2:}  We now prove the general statement.  Let
$(f_1, f_1^{-1}X_1), (f_2, f_2^{-1}X_2) $ be marked covers and
$\phi_0,\phi_1$ as above give an equivalence carrying the arc
$\alpha^1$ to $\alpha^2$.  By replacing $f_2$ with the conjugate
$\phi_0 f_2 \phi^{-1} _0$ we may assume that $\phi_0 = \id$ and that
$f_1(\alpha^1) = f_2(\alpha^2)$.  

By Step 1, the combinatorial class of $g_i$ is independent of the
choice of discs $W, D$.  Let $W_1,D_1$ be any such choice for $g_1$.
Hence we may assume that $(W_2,D_2) = \phi_1 (W_1,D_1)$ and that
$f_1(\cl{W}_1, \alpha^1) = f_2(\cl{W}_2,\alpha^2)$.  

Let $\psi_0: (S^2,X)) \to (S^2,X)$ be the identity map (after
conjugation, we may assume $X=X_1 = X_2$ ).  We will find a
homeomorphism $\psi_1$ such that $\psi_0 g_1 = g_2 \psi_1$ and
$\psi_0$ is isotopic to $\psi_1$ relative to $X$.  As in the
proof of the previous step, we define $\psi_1$ in four pieces $\psi_1
| S^2 - W_1 = \psi^{(1)}, \psi^{(2)}, \psi^{(3)}, \psi^{(4)}$.  The
details are all similar to those in the proof of the previous step.
\qed

\rmk A similar proof shows that when blowing up along any finite collection
of arcs $\alpha^1, \alpha^2, ...$ of arcs for which $\Int(\alpha^i)
\cap \alpha^j = \vac, i \neq j$, the equivalence class of the
resulting map is independent of the choices and 
representatives.  Moreover, the resulting class is the same whether we
blow up one arc at a time or all of them simultaneously.  

The next proposition gives a condition under which arc systems in
$(S^2,X)$ forward invariant under $f$ up to isotopy are also forward
invariant under $g$ up to isotopy.  

\begin{prop}[Orbits off $\alpha$ persist]
\label{prop:orbits off alpha persist}
Suppose $\Lambda$ is an arc system in $(S^2,X)$ which is forward
invariant under $f$ up to isotopy relative to $X$.  If, under $f$,
there exists a lift $\wt{\Lambda} \simeq_X \Lambda$ such that
$\wt{\Lambda} \cap \Int(\alpha) = \vac$, then $\Lambda$ is forward
invariant under $g$ up to isotopy relative to $X$. 
\end{prop}

\pf Let $\Lambda_0 \subset \Lambda$ be an arc system and
$\wt{\Lambda}_0$ a lift of $\Lambda_0$ under $f$ such that
$\wt{\Lambda}_0 \simeq_X \Lambda$.  Since $g|\rs - D = f \circ c$, if
$\wt{\lambda}_0$ is a lift under $f$ of $\lambda_0 \in \Lambda_0$ such
that $\wt{\lambda}_0 \cap \Int(\alpha) = \vac$, then $h_1
(\wt{\lambda}_0) = c^{-1}(\wt{\lambda}_0)$ is a lift under $g$ of
$\lambda_0$.  Since $h_1$ is isotopic to the identity relative to $X$,
$\wt{\lambda}_0 \simeq_X h_1 (\wt{\lambda}_0)$.  Hence $h_1
(\wt{\Lambda}_0)$ is a lift of $\Lambda_0$ under $g$ such that $h_1
(\wt{\Lambda}_0) \simeq_X \Lambda_0$.  
\qed

\section{Proofs of Theorems A and B}
\label{section:proofs}

In \S 3.1 we state the characterization of marked rational maps as
marked branched coverings as an immediate corollary to the proof of
Thurston's characterization of rational maps as branched coverings.
In \S 3.2 we state the main ingredient in our proof, Theorem
\ref{thm:arcs intersecting obstructions}, {\it Arcs intersecting
obstructions}; its proof we relegate to an Appendix. In
\S 3.3 and we prove Theorem A; \S 3.4 gives the proof of Theorem B, 
which depends on a construction given in \S 3.5.

\subsection{Characterization of marked rational maps}

Suppose
$(f,X)$ is a marked branched covering.  To a curve system $\Gamma$
in $(S^2, X)$ we associate the {\it Thurston linear transformation}
which is the linear map
\[	f_{\Gamma}: \R ^{\Gamma} \rightarrow \R ^{\Gamma}	\]
given by 
\[
f_{\Gamma}(\gamma) =  \sum_{\gamma' \subset f^{-1}(\gamma)}
\frac{1}{\deg (f: \gamma' \to \gamma) }[\gamma']_{\Gamma},
\]
where $[\gamma']_{\Gamma}$ denotes the element of $\Gamma$ isotopic to
$\gamma'$, if it exists.  If there are no such elements, the sum is
taken to be zero.  Here, we regard elements of $\Gamma$ as the basis
vectors of $\R^{\Gamma}$.  This transformation depends only on the
isotopy classes of curves relative to $X$ and is natural with respect
to iteration of $f$.  Since the entries of $f_{\Gamma} $ are
nonnegative, the Perron-Frobenius Theorem (\cite{gantmacher}, Chapter
XIII) implies that there is a nonnegative real eigenvalue
$\lambda(f_{\Gamma})$ equal to the spectral radius of $f_{\Gamma}$.  

A nonnegative square matrix $A_{ij}$ is said to be {\it irreducible}
if for each $(i,j)$, there is an $n \geq 0$ such that $A^n _{ij} >0$.
We say that $\Gamma$ is irreducible if the matrix for $f_{\Gamma}$ is
irreducible.  The Perron-Frobenius theory easily gives that if
$\Gamma$ is any curve system with $\lambda(f_\Gamma) >0$, then then
$\Gamma$ contains an irreducible curve system $\Gamma'$ for which
$\lambda(f_{\Gamma'}) = \lambda(f_{\Gamma})$.  

A curve system $\Gamma$ is said to be a {\it multicurve} if its
elements are nonperipheral and pairwise disjoint.  The proof of
Thurston's characterization of rational maps as branched coverings of
the sphere given in \cite{DH1} extends immediately to the setting of
branched covers marked by additional points.

\begin{thm}[Characterization of marked rational maps]
\label{thm:characterization of marked rational maps}
Let $(f,X)$ be a marked branched covering.  Then $(f,X)$ is
combinatorially equivalent to a marked rational map $(R,Y)$ if and only
if for every irreducible multicurve $\Gamma$ in $(S^2,X)$, either

\begin{enumerate}

\item $\lambda(f_{\Gamma})=1$ and the orbifold associated to $(f,P(f))$ has
signature $(2,2,2,2)$, in which case $R$ is covered by an integral
endomorphism of the torus, or 

\item $\lambda(f_{\Gamma}) < 1$, in which case the rational map $R$ is unique up
to conjugation by elements of $\Aut ( \rs )$.

\end{enumerate}
\end{thm}

\begin{defn}[Thurston obstruction.]  
If $(f,X)$ is a marked branched covering for which the orbifold
associated to $(f,P(f))$ is not the $(2,2,2,2)$ orbifold, an
(irreducible) multicurve $\Gamma$ for which $\lambda(f_{\Gamma}) \geq
1$ is called an (irreducible) {\it Thurston obstruction}.  
\end{defn}

\noindent {\bf Example:} Let $\Gamma = \{\gamma_0, ...,\gamma_{p-1}\}$ be a
multicurve such that the following holds: for each $0 \leq j \leq
p-1$, $\gamma_{j}$ has a (not necessarily unique) lift $\delta
\simeq_X \gamma_{j-1 (\bmod p)}$ and $\deg{f:\delta
\to \gamma_j} = 1$.  Then $\Gamma$ is called a {\it Levy cycle} and 
is an example of an irreducible Thurston obstruction.   

\gap\noindent{\bf Remarks:} 
\begin{enumerate}
\item See \cite{DH1} for the proof and for the definition of the orbifold
associated to a rational map. A rational map $R$ is said to be {\it
covered by a torus endomorphism} if there is a complex torus $T$, a
holomorphic covering map $\wt{R}: T \to T$, and a surjective
holomorphic map $p$ such that $R \circ p = p \circ \wt{R}$.  The
endomorphism $\wt{R}$ is said to be {\it integral} if it has a lift to
an automorphism of $\C$ of the form $z \mapsto n z$, some $n \geq 2$.
The first case in the above theorem is rare and does not arise, for
example, if there are periodic critical points of $f$ or if $|P(f)|
\geq 5$.

\item In \cite{DH1} the curve systems $\Gamma$ in the theorem are 
required to be {\it $f$-invariant}, i.e. every lift $\delta$ of an element
$\gamma \in \Gamma$ is either inessential, peripheral, or isotopic to
an element of $\Gamma$.  However, by taking preimages we may always
extend an irreducible multicurve $\Gamma$ to an $f$-invariant
multicurve without decreasing the leading eigenvalue of the
corresponding Thurston transformation.
\end{enumerate}

To apply Theorem \ref{thm:characterization of marked rational maps} we
will need to rule out the first rare case.  Let $(f,X)$ be a marked
rational map $f(z)$.  If $f$ is covered by an integral torus
endomorphism, then $f$ is not a polynomial.  Moreover, there can be no
arcs in $(S^2,X)$ with distinct endpoints which are periodic under $f$
up to isotopy relative to $X$.  For the Euclidean metric on $\C$
pushes down to an (orbifold) metric on $\rs$ which is uniformly
expanded under $f$.  Hence all inverse branches of $f$ are uniformly
contracting with respect to this orbifold metric.  This implies that
if $\wt{\lambda}$ is any lift under $f$ of an arc $\lambda$ in
$(S^2,X)$, then the diameter of $\wt{\lambda}$ is strictly less than
that of $\lambda$.  Since the diameter of any arc in $(S^2,X)$ with
distinct endpoints is positive, this shows that there are no arcs with
distinct endpoints which are periodic under $f$ up to isotopy relative
to $X$.  Hence to prove our two theorems, we may assume that the
orbifold $O_{f}$ is not the $(2,2,2,2)$ orbifold.

If $g$ is $f$ blown up along $\alpha$ then the Euler characteristic
$\chi(O_{g}) < \chi(O_{f}) \leq 0$.  The $(2,2,2,2)$ orbifold has
Euler characteristic zero, hence $O_{g}$ cannot be the $(2,2,2,2)$
orbifold. So $g$ is combinatorially equivalent to a rational map if
and only if there are no irreducible Thurston obstructions for $g$.

\subsection{Arcs intersecting obstructions}
\label{subsection:arcs intersecting obstructions}

Let $(f,X)$ be a marked branched covering, $\Lambda$ be an arc
system in $(S^{2},X)$, and $\Gamma$ a curve system in $(S^2,X)$.  

The {\it unweighted Thurston linear transformation} $f_{\#, \Gamma}:
\R ^{\Gamma} \to \R ^{\Gamma}$ is defined by
\[ f_{\# , \Gamma}(\gamma) = \sum_{\gamma' \subset f^{-1}(\gamma)}
[\gamma']_{\Gamma} \; \mbox{for $\gamma \in \Gamma$,} \] where
$[\gamma']_{\Gamma}$ denotes the element in $\Gamma$ homotopic to
$\gamma'$ and zero otherwise.  As before, this map is independent of
choice of representatives and natural with respect to iteration.  Note
that $0 \leq (f_\Gamma)_{i,j} \leq (f_{\#,\Gamma})_{i,j}$ for all
$i,j$ and that $(f_\Gamma)_{i,j} = 0$ if and only if
$(f_{\#,\Gamma})_{i,j} = 0$.  This implies that $f_\Gamma$ is
irreducible if and only if $f_{\#,\Gamma}$ is irreducible.

In analogy with the unweighted Thurston linear transformation, we
define
\[ f_{\#, \Lambda}: \R ^{\Lambda} \to \R^{\Lambda} \]
by setting
\[ f_{\#, \Lambda}(\lambda) = \sum_{\lambda' \subset f^{-1}(\lambda)}
[\lambda']_{\Lambda} \; \mbox{for $\lambda \in \Lambda$}, \] where
$[\lambda']_{\Lambda}$ denotes the element of $\Lambda$ homotopic to
$\lambda'$ (if it exists).  We say that $\Lambda$ is {\it irreducible}
if $f_{\# ,\Lambda}$ is irreducible.  It is straightforward to verify
that the transformation $f_{\#,\Lambda}$ depends only on $[\Lambda]_X$
and not on the choice of representatives and that $f_{\#,\Lambda}$ is
natural with respect to iteration of $f$, i.e. $(f^{n})_{\#,\Lambda} =
(f_{\#,\Lambda})^n$.  

If $\alpha$ and $\beta$ are arcs in $(S^2,X)$ we define the 
{\it multiplicity} $\mult(f: \alpha \to \beta)$ to be the number of
lifts of $\beta$ which are isotopic to $\alpha$ relative
to $X$.  Thus $(f_{\#,\Lambda})_{ij}=\mult(f:\lambda_i \to \lambda_j)$.

The following proposition follows immediately from the decomposition
theorem for nonnegative square matrices in \cite{gantmacher}, Chapter
XIII.

\begin{prop}[{\bf Decomposition of invariant arc systems}]
\label{prop:decomposition of invariant arc systems}
Let $\Lambda$ be a finite arc system forward invariant under $f$ up to
isotopy relative to $X$.  Then there are disjoint (as subsets of
$\Lambda$) irreducible arc subsystems $\Theta_1, \Theta_2, ..., \Theta_s$ of
$\Lambda$ and integers $n_1, n_2, ..., n_s$ such that for each
$\lambda \in \Lambda$, there is a (not necessarily unique)
$\Theta_k$ such that $\lambda$ is isotopic to a lift of an element of
$\Theta_k$ under $f^{n_k}$.  
\end{prop}

For $\Lambda$ a curve system or an arc system in $(S^2,X)$, set
$\widetilde{\Lambda}$ (resp. $\widetilde{\Lambda}(f^n)$) to be the
union of those components of $f^{-1}(\Lambda)$
(resp. $f^{-n}(\Lambda)$) which are isotopic to elements of
$\Lambda$. If $\Lambda$ is irreducible, each component of $\Lambda$ is
isotopic to some (not necessarily unique) component of
$\widetilde{\Lambda}(f)$ (resp.  $\widetilde{\Lambda}(f^n)$).

The following theorem appears in a different form in
\cite{shishikura:tan:cubicmating}.  

\begin{thm}[Arcs intersecting obstructions] 
\label{thm:arcs intersecting obstructions}
Let $(f,X)$ be a marked branched covering, $\Gamma$ an irreducible
Thurston obstruction in $(S^{2},X)$, and $\Lambda$ an irreducible arc
system in $(S^{2},X)$. Suppose furthermore that
$\#(\Gamma\cap\Lambda)=\Gamma \cdot \Lambda$.  Then either
\begin{enumerate}

\item $\Gamma \cdot \Lambda = 0$ and $\Gamma\cdot f^{-n}(\Lambda) =
0$, $f^{-n}(\Gamma)\cdot \Lambda=0$ for all $n\ge 1$; or

\item $\Gamma \cdot \Lambda \neq 0$ and
\begin{enumerate}

\item Each component of $\Gamma$ (resp. $\Lambda$) is isotopic to a
unique component of $\widetilde{\Gamma}$
(resp. $\widetilde{\Lambda}$), the mapping $f:\widetilde{\Gamma}\to
\Gamma$ (resp. $f: \widetilde{\Lambda}\to \Lambda$) is a
homeomorphism, $\widetilde{\Gamma}\cap
(f^{-1}(\Lambda)-\widetilde{\Lambda}))=\emptyset$, and
$\widetilde{\Lambda}\cap
(f^{-1}(\Gamma)-\widetilde{\Gamma})=\emptyset$. More precisely,

\begin{enumerate} 

\item for each $\gamma \in \Gamma$, there is exactly one curve
$\gamma' \subset f^{-1}(\gamma)$ such that $\gamma' \cap
\widetilde{\Lambda} \neq \emptyset$.  Moreover, the curve $\gamma'$ is
the unique component of $f^{-1}(\gamma)$ which is isotopic to an
element of $\Gamma$;

\item for each $\lambda \in \Lambda$, there is exactly one connected
component $\lambda'$ of $f^{-1}(\lambda)$ such that $\lambda' \cap
\widetilde{\Gamma} \neq \emptyset$.  Moreover, the arc $\lambda'$ is
the unique component of $f^{-1}(\lambda)$ which is isotopic to an
element of $\Lambda$;

\end{enumerate}

\item the transformations $f_{\# ,\Gamma}$ and $f_{\# ,\Lambda}$ are
transitive permutations of the basis vectors.  

\item The above results remain true if we replace $f$ by $f^n$, for any
$n \geq 1$ (though transitivity may fail).   

\item  For any arc $\lambda$ in $(S^2,X)$ which is isotopic to an
arc in $f^{-n}(\Lambda)$ for some $n>1$ but which is not isotopic to
an arc in $\Lambda$, $\Gamma \cdot \lambda=0$.  Similarly, for any 
curve $\gamma$ in $(S^2,X)$ which is isotopic to a curve in $f^{-n}(\Gamma)$ 
for some $n>1$ but which is not isotopic to a curve in $\Gamma$, $\Lambda \cdot
\gamma = 0$.  
\end{enumerate} 

\end{enumerate} 
\end{thm}

We prove this theorem in the Appendix.

\begin{cor}
If $\Lambda$ contains more than one periodic cycle, or if the
multiplicity of $f$ from some element $\alpha \in \Lambda$ to $\beta
\in \Lambda$ is larger than one, then there are no irreducible
Thurston obstructions which intersect $\Lambda$.
\end{cor}

\pf If so then $f_{\#, \Lambda}$ is not a transitive permutation of
the basis vectors.
\qed

\subsection{Proof of Theorem A}
\label{subsection:proof of theorem one}
By Theorem \ref{thm:characterization of marked rational maps}, {\it
Characterization of marked rational maps}, and our discussion at the
end of Section 3.1, $(g,X)$ is combinatorially equivalent to a marked
rational map if there exist no irreducible Thurston obstructions
$\Gamma$ for $g$.  Let $\Gamma$ be such an obstruction and let
$\widetilde{\Gamma}$ be the union of those components of
$g^{-1}(\Gamma)$ isotopic to elements of $\Gamma$, relative to $X$.

By the second hypothesis of the theorem and Proposition
\ref{prop:orbits off alpha persist}, {\it Orbits off $\alpha$
persist}, for each $j$, $\Lambda_j$ is forward invariant under $g$ up
to isotopy relative to $X$.  

By Proposition \ref{prop:decomposition of invariant arc systems}, {\it
Decomposition of invariant arc systems}, for each $j$, there is an
irreducible arc system $\Theta_j \subset \Lambda_j$ and an integer
$n_j$ such that $f(\alpha_j)$ is isotopic relative to $X$ to a lift of
an element of $\Theta_j$ under $g^{n_j}$.

We now show that we may choose representatives of $\Gamma$ so that
$\widetilde{\Gamma} \cap D = \vac$.  This will complete the proof of
the theorem: for since $g|\rs - (D-e(\alpha)) = f \circ c$ and $c: \rs
- (D-e(\alpha)) \to \rs - Int(\alpha)$ is isotopic to the identity
relative to $X$, $\Gamma$ is then an irreducible Thurston obstruction
for the marked rational map $(f,X)$.  By Theorem
\ref{thm:characterization of marked rational maps}, {\it
Characterization of marked rational maps}, and the fact that $f$ is
not covered by an integral torus endomorphism, this gives a
contradiction.

First, we may assume $\Gamma$ minimizes the number of intersections
with $f(\alpha_j)$ and $\Theta_j$ for all $j$.  Write $\bdry D =
\cup_{j}\alpha_{0,j}^{\pm}$, where $\alpha_{0,j}^{\pm}$ are the two
unique components of $(g|D)^{-1}(f(\alpha_j))$ which are contained in
$\bdry_{\pm} D$ respectively.  The proof now breaks down into several
cases:

\gap
\noindent {\bf Case 1: $\Gamma \cdot f(\alpha_j) = 0$}.  Then $\Gamma
\cap f(\alpha_j) = \vac$, hence $\widetilde{\Gamma} \cap
g^{-1}(f(\alpha_j)) = \vac$, and so $\widetilde{\Gamma} \cap
\alpha_{0,j}^{\pm} = \vac$.

\gap
\noindent {\bf Case 2: $\Gamma \cdot f(\alpha_j) \neq 0$}.  
We will again show that $\wt{\Gamma} \cap \alpha^{\pm}_j = \vac$.  

We first show that in this case $f(\alpha_j) \in \Theta_j$ up to isotopy.  
By our choice of $\Theta_j$, $f(\alpha_j)$ is isotopic to a lift of an
element $\theta \in \Theta_j$ under $f^{n_j}$.  First, if $\Gamma \cdot
\Theta_j = 0$, then by an application of Theorem \ref{thm:arcs
intersecting obstructions}, {\it Arcs intersecting obstructions}, Part
1, we would have $\Gamma \cdot f(\alpha_j) = 0$, contradicting our
assumption.  Hence $\Gamma \cdot \Theta_j \neq 0$.  Next, if
$f(\alpha_j) \not\in \Theta_j$ up to isotopy, then by Theorem
\ref{thm:arcs intersecting obstructions}, {\it Arcs intersecting
obstructions}, Part 2(d) ( with $\Lambda = \Theta_j$ and
$\lambda=f(\alpha_j)$), $\Gamma \cdot f(\alpha_j) = 0$, contradicting
our assumption.  Hence $f(\alpha_j) \in \Theta_j$ up to isotopy.

Neither of the arcs $\alpha_{0,j}^{+}, \alpha_{0,j}^{-}$ can be
isotopic to elements of $\Theta_{j}$.  For if so, then $\alpha$ must
consist of a single arc in $(S^2,X)$, i.e. $L=1$ in the definition of
the blown up map $g$. Then $\alpha$ itself is isotopic to an element of
$\Theta_j$ since $\Theta_j$ is irreducible.  It then follows that
$\bdry D \subset g^{-1}(f(\alpha))$ consists of two arcs
$\alpha_{0}^{+}$ and $\alpha_{0}^{-}$ in $(S^2,X)$ which are isotopic
relative to $X$, by construction of the blown-up map $g$.  Hence the
multiplicity of $g: \alpha_{0}^{\pm} \to f(\alpha)$ is at least $n
\geq 2$, where $n$ is the number of times we blow up $\alpha$.  Hence
the transformation $g_{\#,\Theta}$ is not a transitive permutation of
the basis vectors.  By Theorem \ref{thm:arcs intersecting
obstructions}, {\it Arcs intersecting obstructions}, $\Gamma \cdot
\Theta = 0$, hence $\Gamma \cdot f(\alpha) = 0$, contradicting our 
assumption. 

By the previous paragraph $\alpha_{0,j}^{\pm} \not\in \Theta_{j}$ up
to isotopy.  By Theorem \ref{thm:arcs intersecting obstructions}, Part
2a(ii), there is exactly one lift $\lambda'$ of $f(\alpha_j)$ which
intersects (as a subset of $S^2$) the set $\widetilde{\Gamma}$, and
this lift $\lambda'$ is isotopic to an element of $\Theta_{j}$.  Hence
$\alpha_{0,j}^{\pm} \cap \widetilde{\Gamma} = \vac$.

\gap In both cases we have shown that for all $j$, $\widetilde{\Gamma}
\cap \alpha_{0,j}^{\pm} = \vac$.  Thus $\widetilde{\Gamma} \cap \bdry
D = \vac$,  which implies that $\widetilde{\Gamma} \cap D = \vac$:
there can be no components of $\Gamma$ contained entirely in $D$ since
$\Int(D) \cap X = \vac$.
\qed 

\gap\noindent{\bf Proof of Corollary 1}
Suppose $\alpha$ is a periodic or preperiodic internal ray of $f$.
Points in the interior of an internal ray lie in Fatou components and
have infinite forward orbits under $f$.  Hence no points in the
interior of an internal ray can be preperiodic under $f$.  We may
therefore assume that $\alpha$ consists of a single arc in
$(S^2,f^{-1}X)$, and hence that $f(\alpha)$ is a periodic or
preperiodic internal ray with endpoints in $X$.

Two internal rays either coincide, or intersect only in their
endpoints.  Hence if $l_i = f^i (f(\alpha)), i > 0$, then $l_i \cap 
\Int(f(\alpha)) \neq \vac $ if and only if $\alpha=l_i$ for some $l_i$,
i.e. $f(\alpha)$ is a periodic internal ray. 

Let $\cal{L}$ $ = \{l_i\}_{i>0}$.  Then $\cal{L}$ is a finite collection of
arcs in $(S^2,X)$ which, as a subset of the sphere, is forward
invariant under $f$.  By Proposition \ref{prop:finding invariant arc
systems} there exists a subset $\Lambda \subset \cal{L}$ which is an
arc system forward invariant under $f$ such that each element of
$\cal{L}$ is isotopic relative to $X$ to a unique element of
$\cal{L}$.  By the definition of forward invariance there is a subset
$\Lambda_0 \subset \Lambda$ and a lift $\wt{\Lambda}_0
\simeq_X \Lambda$.   

If $f(\alpha) \not\in \Lambda$, then $\Int(f(\alpha)) \cap \Lambda =
\vac$, hence $\Int(\alpha) \cap f^{-1}(\Lambda) = \vac$.  Since
$\wt{\Lambda}_0 \subset f^{-1}(\Lambda)$, $\wt{\Lambda}_0 \cap
\Int(\alpha) = \vac$.  The Corollary then follows by Theorem A.    

If $f(\alpha) \in \Lambda$, let $\lambda'$ be any arc which is
isotopic to $f(\alpha)$ relative to $X$ and such that $\Int(f(\alpha))
\cap \lambda' = \vac$.  Let $\Lambda'$ be the arc system $(\Lambda -
\{f(\alpha)\}) \cup \{\lambda'\}$.  Then $\Lambda \simeq_X \Lambda'$ and 
so $\Lambda'$ is forward invariant under $f$ up to isotopy relative to
$X$.  By construction $\Lambda' \cap \Int(f(\alpha)) = \vac$, hence
$f^{-1}(\Lambda') \cap \Int(\alpha) = \vac$.  If $\Lambda' _0 $ is the
unique subset of $\Lambda'$ for which $\Lambda' _0 \simeq_X
\Lambda_0$, and if $\wt{\Lambda}' _0$ is the lift of $\Lambda' _0$
obtained by lifting isotopies between elements of $\Lambda' _0$ and
$\Lambda_0$, then $\wt{\Lambda}' _0 \simeq_X \wt{\Lambda}_0 \simeq_X
\Lambda$.  Since $\wt{\Lambda}'_0 \subset f^{-1}(\Lambda')$,
$\wt{\Lambda}' _0 \cap \Int(\alpha) = \vac$.  The Corollary then
follows by Theorem A.

\qed

\noindent{\bf Proof of Corollary 2}  The argument is similar to the one
given above.  The hypothesis guarantees the existence of the arc $\lambda'$
used in the last paragraph.
\qed

\subsection{Proof of Theorem B}
\label{subsection:proof of theorem two}

Let $f(z)$ be a postcritically finite monic polynomial and $X$ a
finite set such that $f(X) \subset X$ and $X \supset P(f)$ .  Let
$K(f)$ and $F(f)$ denote the filled-in Julia set and Fatou set of $f$,
respectively.

In \cite{poirier:pcfpii}, \S I.1 (originally in \cite{DH2}), it is
shown that $X-\{\infty\}$ determines uniquely a finite topological
tree $T$ with vertices $V(T)$ satisfying the following properties.

\begin{enumerate}
\item $T \subset K(f)$ and $T \cap F(f)$ only in internal rays;
\item $f(T) \subset T$ and $f(V(T)) \subset V(T)$;
\item $V(T) \supset X-\{\infty\}$.
\end{enumerate}

The definition of internal and external rays is given in the following
section.  

Our first step in the proof is to construct a ``spider'' $S$ for $T$.
The construction of $S$ is given in the following section.  $S$ is a
finite topological graph with vertices $V(S)$ and edges $E(S)$
satisfying the following properties:

\begin{enumerate}
\item $V(S) = X \supset P(f)$;
\item Every edge of $S$ is a topological arc joining a finite vertex
to infinity; 
\item $f(S) \subset S$;
\item $f$ maps edges of $S$ homeomorphically to edges of $S$;
\item $S \cap T = X-\{\infty\}$;
\end{enumerate}

Note that the edges of $S$ are arcs in $(S^{2},X)$.

Given the existence of the spider $S$ with the properties above, we
now prove Theorem B.  

The edges $E(S)$ form a finite collection of arcs $l$ in $(S^2,X)$
such that $f(l) \in E(S)$ for each $l \in E(S)$.  By Proposition
\ref{prop:finding invariant arc systems}, there is an arc system $\Lambda
\subset E(S)$ which is forward invariant under $f$ up to isotopy
relative to $X$, i.e. there is a subset $\Lambda_0 \subset \Lambda$
and a lift $\wt{\Lambda}_0 \simeq_X \Lambda$.  Note that since
$\wt{\Lambda}_0 \simeq_X \Lambda \subset X$, the endpoints of
$\wt{\Lambda}$ are points in $X$. 

We now show that $\Lambda$ is forward invariant under $g$ up to
isotopy relative to $X$.  By Proposition \ref{prop:orbits off alpha
persist}, {\it Orbits off $\alpha$ persist}, it is enough to show that
$\wt{\Lambda}_0 \cap \Int(\alpha) = \vac$.  Since $\alpha \subset T$
and $f(T) \subset T$, $f(\alpha)\subset T$.  Since $T \cap S \subset
X$, $f(\alpha) \cap S \subset X$.  Hence $\alpha \cap f^{-1}(S) 
\subset f^{-1}(X)$.  Now, $\wt{\Lambda}_0
\subset f^{-1}(S)$, and as remarked above, the endpoints of
$\wt{\Lambda}_0$ are points in $X$, hence $\alpha \cap \wt{\Lambda}_0 
\subset X$.  But $\Int(\alpha) \cap X = \vac$, hence $\Int(\alpha) 
\cap \wt{\Lambda}_0 = \vac$. 

Hence $\Lambda$ is forward invariant under $g$ up to isotopy relative
to $X$. By Proposition \ref{prop:decomposition of invariant arc
systems}, there are disjoint (as subsets of $\Lambda$) irreducible arc
systems $\Theta_k \subset \Lambda$ and integers $n_k$ such that each
element of $\Lambda$ is a lift of an element of $\Theta_k$ under
$n_k$ for some $k$.

Let $\Gamma$ be an irreducible Thurston obstruction for $g$ and let
$\widetilde{\Gamma}$ be the union of components of $g^{-1}(\Gamma)$
which are isotopic relative to $X$ to elements of $\Gamma$.  
We may assume that $\Gamma$ minimizes the number of intersections with
each $f(\alpha_j)$ and with each $\lambda \in \Lambda$.  As in the
previous proof it suffices to show that $\widetilde{\Gamma}
\cap \bdry D = \vac$, where $D$ is the disc used in the construction
of $g$.  

Let $\cal{K}$ $ = (g|D)^{-1}(\Lambda)$; this is a finite, connected graph
in $D$ with vertices equal to $g^{-1}(X) \cap D \subset \cal{K}$.

We first reduce to the case when $\widetilde{\Gamma} \cap \cal{K} =
\vac$.  As in the previous proof, there are two cases.  Let $\lambda
\in \Lambda$. 
\gap

\noindent {\bf Case 1:} $\Gamma \cdot \lambda = 0$.  Then $\Gamma \cap
\lambda = \vac$ by our choice of representatives of $\Gamma$, so
$\widetilde{\Gamma} \cap (g|D)^{-1}(\lambda)= \vac$.
\gap

\noindent {\bf Case 2:} $\Gamma \cdot \lambda \neq 0$.  

Let $\wt{\lambda}$ be any lift of $\lambda$ under $g$ which is
contained in $D$ (and hence in $\cal{K}$).  By Proposition
\ref{prop:decomposition of invariant arc systems}, $\lambda$ is
isotopic to a lift of an element of some irreducible arc system
$\Theta_j \subset \Lambda$ under $f^{n_k}$, some $n_k$. 

Since $\lambda \in \Lambda \subset E(S)$, $\lambda$ is an edge of the
spider $S$ and hence joins a point of $X-\{\infty\}$ to the point
$\infty$.  By construction $D \cap \infty = \vac$.  Since
$\wt{\lambda} \subset D$ is a lift of $\lambda$ under $g$, one of the
endpoints of $\wt{\lambda}$ is contained in $D$ and maps onto $\infty$.
Hence this endpoint is strictly preperiodic under $g$.  Hence
$\wt{\lambda}$ cannot be contained in the irreducible arc system
$\Theta_k$.  By Theorem \ref{thm:arcs intersecting obstructions}, {\it
Arcs intersecting obstructions}, $\wt{\Gamma} \cap \wt{\lambda} = \vac$.  

Together, the two cases show that $\wt{\Gamma} \cap \cal{K}=\vac$. 

\gap

So we may assume $\widetilde{\Gamma} \cap \cal K = \vac$.  We conclude
the proof by showing that if $\wt{\Gamma} \cap \bdry D \neq \vac$,
then $\Gamma$ does not have minimal intersection number with
$f(\alpha)$, contrary to our assumption.

Suppose $\wt{\gamma} \subset \wt{\Gamma}$ intersects $\bdry D$. Since
$\wt{\Gamma} \cap \cal{K} = \vac$, $\wt{\gamma} \cap \Int(D)$ is
contained in a component $V$ of $D-\cal{K}$.  Since $V(S) = X \supset
P(f)$, and since each edge $E(S)$ is isotopic relative to $X$ to an
element of some $\Lambda_k$, each component of $S^2 - \Lambda$ is an
open disc containing no critical values of $g$.  Hence each component
of $g^{-1}( S^2 - \Lambda)$ maps homeomorphically under $g$ to its
image.  Since $V$ is contained in such a component, $g|V$ is a
homeomorphism.  Let $W \subset V$ be any component of $V -
\wt{\gamma}$ which ``outermost'', i.e. not separated from $\bdry D$ by
any component of $\wt{\gamma} \cap D$.  Then $W$ is a disc in
$S^2-g^{-1}X$ bounded by exactly one subarc of $\wt{\gamma}$ and one
subarc of $\bdry_+ D$ or $\bdry_- D$.  Since $g|V$ is a homeomorphism
and $W \subset V$, $g(W)$ is a disc in $S^2-X$ bounded by one subarc
of $f(\alpha)$ and one subarc of $\gamma \subset \Gamma$,
contradicting the assumption that $\Gamma \cap f(\alpha)$ is
minimized.
\qed

\subsection{Construction of $S$}
\label{subsection:construction of S}

We first recall some facts about the dynamics of postcritically finite
polynomials from \cite{MIL1}.  

Let $f(z)$ be a postcritically finite polynomial.  The {\it filled-in
Julia set} $K(f)=\{z | f^{n}(z) \not \to \infty\}$ is a connected,
compact subset of $\C$ whose complement is connected.  The Julia set
$J(f) = \bdry K(f)$ is locally connected.  For each Fatou component
$\Omega$, there is exactly one point $x \in \Omega$ such that
$f^{n}(x) \in P(f)$ for some $n$.  Let $\Omega_{x}$ denote the Fatou
component containing such a point $x$.  A classical theorem of
B\"{o}ttcher implies the following: there are holomorphic isomorphisms
$\phi_{x}: (\Delta,0) \to (\Omega_{x},x)$ such that for all $w \in
\Delta$, $\phi_{f(x)}(w^{d_{x}}) = f(\phi_{x}(w))$, where $d_{x}$ is
the local degree of $f$ near $x$.  Since $J(f)$ is locally connected,
each Fatou component has locally connected boundary, and a theorem of
Carath\'eodory implies that the maps $\phi_{x}$ extend continuously
to $S^{1}=\bdry \Delta$.  Let $R_{t}$ denote the radial line $\{r
\exp(2\pi i t) | 0 \leq r \leq 1\}$.  The image
$R_{x,t}=\phi_{x}(R_{t})$ we call the {\it ray of angle $t$ in
$\Omega_{x}$}; if $x=\infty$ $R_{x, t}$ is called an {\it external
ray}; otherwise we call it an {\it internal ray}.  Each periodic point
$x \in J(f)$ is the landing point of at least one and at most finitely
many external rays.  The landing point of an external ray is periodic
if and only if $x$ is periodic, though the periods may differ.
Internal and external rays map homeomorphically onto their images.

\gap\noindent {\bf Construction of $S$.}  Our spider will be a union 
\[ S = \bigcup_{x \in X-\{\infty\}} E_{x}, \]
where $E_{x}$ denotes the set of edges of the graph $S$ incident to
$x$. 
\begin{itemize}
\item If $x \in J(f) \cap (X-\{\infty\})$, each element $e \in E_{x}$ 
will be an external ray landing at $x$, and $E_{x}$ will be the union
of all external rays landing at $x$.

\item If $x \in F(f) \cap (X-\{\infty\})$, each element $e \in E_{x}$ 
will be the union of exactly one internal ray in $\Omega_{x}$ and one
external ray landing at a common point $q \in \bdry \Omega_{x}$.  The
set $E_{x}$ will be a finite union of such pairs.
\end{itemize}

Since external and internal rays map homeomorphically onto their
images, and the image of an internal (external) ray is again an
internal (external) ray, edges will map homeomorphically to edges.

Each $x \in J(f) \cap (X-\{\infty\})$ is eventually periodic, so the
set $E_{x}$ of external rays landing at $x$ is finite (see
\cite{MIL1}, \S 18), and $f(E_{x}) = E_{f(x)}$.  Since $T \subset K(f)$,
$E_{x} \cap T \subset X-\{\infty\}$.

We now define $E_{x}$ for $x \in F(f) \cap (X-\{\infty\})$.  First,
choose one element $x$ from each periodic cycle in $(X-\{\infty\})
\cap F(f)$.  Let $p$ be the period of $x$.  Then $T \cap \Omega_{x}$
consists of a finite collection of eventually periodic internal rays.
Thus there exists a periodic internal ray $R_{x}$, of some period $kp>0$,
such that $R_{x} \cap T = \{x\}$.  Let $q$ be the landing point of the
internal ray $R_{x}$; then $q$ is periodic of period $k$.  There
exists a periodic external ray $R_{\infty}$ landing at $q$.  Since $q$
is not a critical point, $f^{k}$ is a local homeomorphism near $q$
sending $R_{x}$ to itself, hence the ray $R_{\infty}$ also has period
$k$.  Let $e_{x}=R_{x} \cup R_{\infty}$; this will be one edge in the
collection $E_{x}$.  Define
\[ E_{x}=\bigcup_{n=0}^{k-1} f^{np}e(x). \]
$E_{x}$ is finite since $q$ is periodic.  For $y \in f^{n}(x), 1 \leq n
< p$, define
\[ E_{y} = f^{n}(E_{x}). \]
Thus $E_{y}$ is finite.  Since $R_{x}$ is periodic, $R_{x} \cap T =
\{x\}$, $f(X-\{\infty\}) \subset X-\{\infty\}$, and $f(T) \subset T$,
we have that $E_{y} \cap T \subset X-\{\infty\}$ for all $y$.

We now inductively define $E_{x}$ for strictly preperiodic $x \in X-\{\infty\}
\cap F(f)$.  Suppose the collection $E_{y}$ has been defined, where
$f(x)=y$.  Choose any element $e_{y} \in E_{y}$.  The edge $e_{y}$ is a union
of an internal ray $R_{y}$ and external ray $R_{\infty}$; let $q$
be their common landing point.  Choose any preimage $R_{x}$ of
$R_{y}$ joining $x$ to a point $q'$ in $\bdry \Omega_{x}$.  Since a
small disc near $q$ pulled back to $q'$ is again a disc, there is a
preimage $R'_{\infty}$ of $R_{\infty}$ joining $\infty$ to $q'$.
Let $e_{x} = R_{x} \cup R'_{\infty}$ and $E_{x}=\{e_{x}\}$.  Then $f:
e_{x} \to e_{y}$ and $E_{x} \cap T = \{x\}$.

This completes the construction of $S$.

\section{Complements}
\label{section:complements}

\subsection{Blowing up multiple arcs}

Theorems A and B generalize in the obvious way to blowing up along
finite collections $\alpha^1, \alpha^2, ...$ of arcs for which $\Int(\alpha^i) 
\cap \alpha^j = \vac, i \neq j$.  

\subsection{Blowing up M\"{o}bius transformations}
\label{subsection:blowing up moebius transformations}
Let $G$ be a finite connected graph in $S^2$ with vertex set $X$ and
edge set $E$.  Suppose in addition that no edge joins a vertex to
itself and that no two vertices are joined by more than one edge.
Assign to each edge $e \in E$ an integer $m(e) \geq 0$.  Let $(M,X)$
be a marked M\"{o}bius transformation of finite order and suppose
$M(G)=G$.

\gap\noindent{\bf Theorem: (Blowing up M\"obius transformations)}
\label{thm:blowing up moebius transformations}
{\it Let $(g,X)$ be the marked branched covering $(M,X)$ blown up $m(e)$
times along the edge $e \in E$.  Let $G'$ denote the subgraph of $G$
which is the union of edges $e'$ (and their vertices) for which
$m(M^{i}e')>0$ for some $i$.  Then $(g,X)$ is equivalent to a marked
rational map if and only if 
\begin{enumerate}
\item $G'$ is connected, and
\item if $G'$ is nonempty, no component of $S^{2}-G'$ contains more 
than one point of $X$.
\end{enumerate}
}

\pf If $G'$ is empty there is nothing to prove.  Otherwise, there is
at least one edge which is blown up.  The vertices of this blown up
edge become periodic critical points for $g$, hence the orbifold of
the (unmarked) branched covering $g$ is not the $(2,2,2,2)$ orbifold.
Hence by Theorem \ref{thm:characterization of marked rational maps},
{\it Characterization of marked rational maps}, $(g,X)$ is equivalent
to a marked rational map if and only if $(g,X)$ has no Thurston
obstructions.  

Note that $G'$ is invariant under $M$.  Let $X'$ denote the vertices
of $G'$.  Then for each connected component $G''$ of $G'$, $G''$ has
at least two vertices.  

We now prove the necessity of (1) and (2). First, assume that $(g,X)$
is equivalent to a marked rational map.  Then $(g^{\circ n}, X)$ is
also equivalent to a marked rational map for every $n \geq 1$.  

If (2) fails, then there is a component $U$ of $S^2-G'$ such that $|U
\cap (X - X')| \geq 2$.  Let $\gamma$ be any simple closed curve in
$U$ separating $G'$ from $U \cap (X - X')$.  Then $\gamma$ is
nonperipheral since $|X'| \geq 2$ and $|U \cap(X-X')| \geq 2$ by
assumption.  Choose $n \geq 1$ so that $M^{\circ n} (\gamma) =
\gamma$.  We may further assume that $\gamma$ is disjoint from the
discs $W$ within which we alter the map $(M,X)$ to obtain the map
$(g,X)$.  It then follows that $g^{\circ n}(\gamma) = \gamma$ and that
$(g^{\circ n})|\gamma$ is a homeomorphism.  Then $\gamma$ is a Levy
cycle for $(g^{\circ n},X)$, so $(g^{\circ n},X)$ is not equivalent to
a rational map.

If (1) fails, then there is a simple closed curve $\gamma$ in $S^2 -
(G' \cup X)$ separating $G'$ into two subsets, each of which has at
least two elements of $X$.  Hence $\gamma$ is nonperipheral.  The
same argument as in the previous paragraph then shows that $(g^{\circ
n},X)$ is not equivalent to a rational map.

We now prove the sufficiency of (1) and (2).  We may assume $|X| \geq
4$, since otherwise there are no nonperipheral simple closed curves in
$(S^2,X)$ and hence no Thurston obstructions.  Next, we may assume
that for every $e \in E$ and for every $i > 0$, $M^i e \simeq_X e $ if
and only if $M^i e = e$, by our requirement that there is at most one
edge between any two vertices and the assumption $M(G)=G$.  

For each orbit of an edge in $E'$ under $M$, choose one edge $e' _j$.
Let $\epsilon_j \simeq_X e' _j$ be an arc in $(S^2,X)$ chosen so that
$\epsilon_j \cap G = e(\epsilon_j)$.  Then in particular,
$\Int(\epsilon_j) \cap e' = \vac$ for each $e' \in E'$.  Let
$\Lambda_j$ denote the arc system which is the orbit of $\epsilon_j$
under $M$.  Then $\Lambda_j$ is forward-invariant under $M$ and the
interior of every element of $\Lambda_j$ is disjoint from the edges
along which we blow up.  Hence by Proposition 
\ref{prop:orbits off alpha persist}, {\it Orbits off $\alpha$
persist}, each $\Lambda_j$ is forward-invariant up to isotopy under
$(g,X)$.  

The transformations $g_{\#, \Lambda_j}$ are irreducible but are not
transitive permutations of the basis vectors.  The transformation
$M_{\#, \Lambda_j}$ is a transitive permutation of the basis vectors
since $\Lambda_j$ is the orbit of a single arc under $M$ and is
therefore irreducible.  Hence as matrices $g_{\#, \Lambda_j} \geq
M_{\#, \Lambda_j} \geq 0$, so $g_{\#, \Lambda_j}$ is irreducible since
$M_{\#, \Lambda_j}$ is irreducible.  Moreover, if $\lambda \in
\Lambda_j$ and $\lambda \simeq_X e'$ where $m(e') \geq 1$, then
$\mult(g: \lambda \to g(\lambda)) = m(e') + 1 \geq 2$, and so the
entry of $g_{\#, \Lambda_j}$ corresponding to the ordered pair
$(\lambda, g(\lambda))$ is at least two.  

Let $\Gamma$ be an irreducible Thurston obstruction for $(g,X)$.  We
may assume that the elements of $\Gamma$ are chosen so that $\Gamma
\cdot \Lambda_j = | \Gamma \cap \Lambda_j|$ for each $j$.  By Theorem
\ref{thm:arcs intersecting obstructions}, {\it Arcs intersecting
obstructions}, $\Gamma \cap \Lambda_j = \vac$ for each $j$.  Since
$\cup_j \Lambda_j = G'$, $\Gamma \subset S^2 - G'$.  Conditions (1)
and (2) imply that every component of $S^2-(G' \cup X)$ is either
simply-connected or is a once-punctured disc.  Hence there are no
nonperipheral simple closed curves in $S^2-(G' \cup X)$, and so there
are no Thurston obstructions for $(g,X)$.

\qed

\gap\noindent{\bf Example:}  
Let $f_{p/q}(z)=z^{2}+c_{p/q}$, $q \geq 2$, where $p<q$ and $p$ and $q$ are
relatively prime, where $c_{p/q}$ is the center of the hyperbolic
component tangent to the main cardioid in the $p/q$ limb of the
Mandelbrot set.  We call such a polynomial ``starlike''.  (There is a
topological tree $T$ which is a union of $q$ internal rays, each
joining a point in the finite superattracting cycle containing $0$ to
a common repelling fixed point.  Thus $T$ is a ``star''.  The tree $T$
is mapped homeomorphically under $f$ into itself with rotation number
$p/q$.)  Let $M_{p/q}=e^{2\pi i p/q}z$ and $X_{p/q}=\{\infty\} \cup_{n
\in \Z }e^{2 \pi i n p/q}$.  Let $G$ be the orbit of the arc $[1,
\infty]$ under $M_{p/q}$.  Then $f_{p/q}$ is equivalent to
$(M_{p/q},X_{p/q})$ blown up once along $\alpha = [1,\infty] \subset
G$.
\gap

\rmk Let $M_{p/q}$ be as above.  Let $Y_{p/q}=X_{p/q}-\{\infty\}$ and let 
$G$ be the  unit circle.  Let $\alpha$ be the circular segment joining 
$1$ and $e^{2 \pi i \frac{1}{q}}$.  Then blowing up $(M,Y_{p/q})$ once along
$\alpha$ yields a quadratic non-polynomial rational map.  Mary Rees
(personal communication) has informed us that all quadratic non-polynomial
maps with obstructed tunings arise in this fashion.

\subsection{Killing obstructions by blowing up}

The blowing up construction can be applied to marked branched covers
which are not equivalent to rational maps in such a way that the
result is equivalent to a rational map.

Let $p(z)=z^{2}-1$.  Extend $p$ to a degree two branched cover of the
closed disc $D$ by adjoining the map $z \to z^{2}$ along the circle at
infinity $S^{1}_{\infty}$.  Now let $f$ be the branched cover obtained
by gluing the action of $p$ on two copies $D_{+},D_{-}$ of $D$ along
$S^{1}_{\infty}$.  The map $f$ is called the {\it mating} of $f$ with
itself and is not equivalent to a rational map.  For let $R_{t}$
denote the external ray of angle $t$ for $p$.  Then after adjoining on
$S^{1}_{\infty}$, $R_{t}$ becomes a topological arc $R_{t}^{\infty}$
in $D$ joining $\infty \cdot e^{2\pi i t}$ to the landing point of
$R_{t}$ in $J(p)$.  For the map $f$, the unions of arcs
$R_{1/3}^{\infty}$, $R_{2/3}^{\infty}$ in $D_{+}$ and $D_{-}$ are
joined together to form a simple closed curve $\gamma \subset
S^{2}-P(f)$ which maps to itself homeomorphically by degree one, and
hence is an obstruction.  Note that the interval $[-1,0] \subset
D_{-}$ is mapped homeomorphically to itself under $f$.  Consider the
map $g=(f,P(f))$ blown up once along $\alpha = [-1,0] \subset D_{-}$.
Then by Theorem \ref{thm:arcs intersecting obstructions}, {\it Arcs
Intersecting obstructions}, no obstruction can intersect $\alpha$.
Hence any obstruction for $g$ is isotopic relative to $P(g)$ into
$D_{+}$, where there are no such obstructions.  So $g$ is equivalent to
a rational map while $f$ is not.

\subsection{Generalized matings}
\label{subsection:generalized matings}

One may think of the map $g$ in the previous example as the polynomial
$z^2-1$ mated with the rational map which is $z^2-1$ blown up along
$\alpha$.  (The resulting branched covering is well- defined since
there are no critical points in the interior of the disc $D$ in which
we glue new dynamics.)  Using techniques similar to the ones employed
in the proofs of our theorems one may prove the following.  To set up
the statement, suppose $f_1 (z)$ is a starlike quadratic polynomial
with tree $T$.  Let $\alpha$ be the union of two consecutive edges in
$T$.  Then $\alpha$ satisfies the hypothesis of Corollary 2, so
blowing up $n$ times along $\alpha$ yields a branched covering
equivalent to a rational map $r_n (z)$.  The map $r_n (z)$ has a fixed
critical point of local degree two which we may take to be at
infinity.  

\begin{thm}[Matings with blown-up starlike polynomials exist]
Let $f_2 (z)$ be an arbitrary postcritically finite quadratic
polynomial.  Then the mating of $f_2 (z)$ with $r_n (z)$ is equivalent
to a rational map $g(z)$.
\end{thm}

Similar results can be formulated for higher-degree maps as well as
for generalizations of tunings.  The proofs are all similar to those
given here: an obstruction is forced to lie outside the region on the
sphere in which the new dynamics is glued by Theorem \ref{thm:arcs
intersecting obstructions}, {\it Arcs intersecting obstructions}.  

\subsection{Special maps}

The following theorem is useful since there are interesting examples
of maps produced by blowing up arcs to which Theorems A' and B' do not
apply.  

\gap\noindent{\bf Theorem: (Special maps)} {\it If $g$ is $(f,X)$
blown up along $\alpha$, and if $|P(g)|=3$, then
$g$ is equivalent to a rational map.}

\pf The orbifold of $g$ is hyperbolic, so $g$ is equivalent to a
rational map if and only if there are no Thurston obstructions, which
do not exist since there are no nonperipheral curves in $(S^{2},P(G))$.
\qed

\section{Examples}
\label{section:examples}

In this section we give some examples and applications.  

\subsection{Examples illustrating Theorems}
\label{subsection:examples illustrating theorems}

\leavevmode

\vspace{5pt}
\noindent{\bf Blowing up a periodic arc in the ``rabbit''.}  Let
$f(z)=z^{2}+c$, where $c \approx -0.12256117 + 0.74486177i$ is chosen
so that $0$ is periodic of period three and Im($z$)$ > 0$.  Let
$X=\{0,c,f(c)\}$.  Then $f$ is ``starlike'', so there is a tree
$T$ (shown in Figure
\ref{fig:rabbit_tree.ps}) with vertices consisting of the set $X$ and the
alpha fixed point $p$ of $f$ and which is mapped homeomorphically to
itself under $f$ with rotation number $1/3$.  $T$ is a union of the
three internal rays for $f$ which are fixed under $f^{\circ 3}$. 
 
\yincludeps{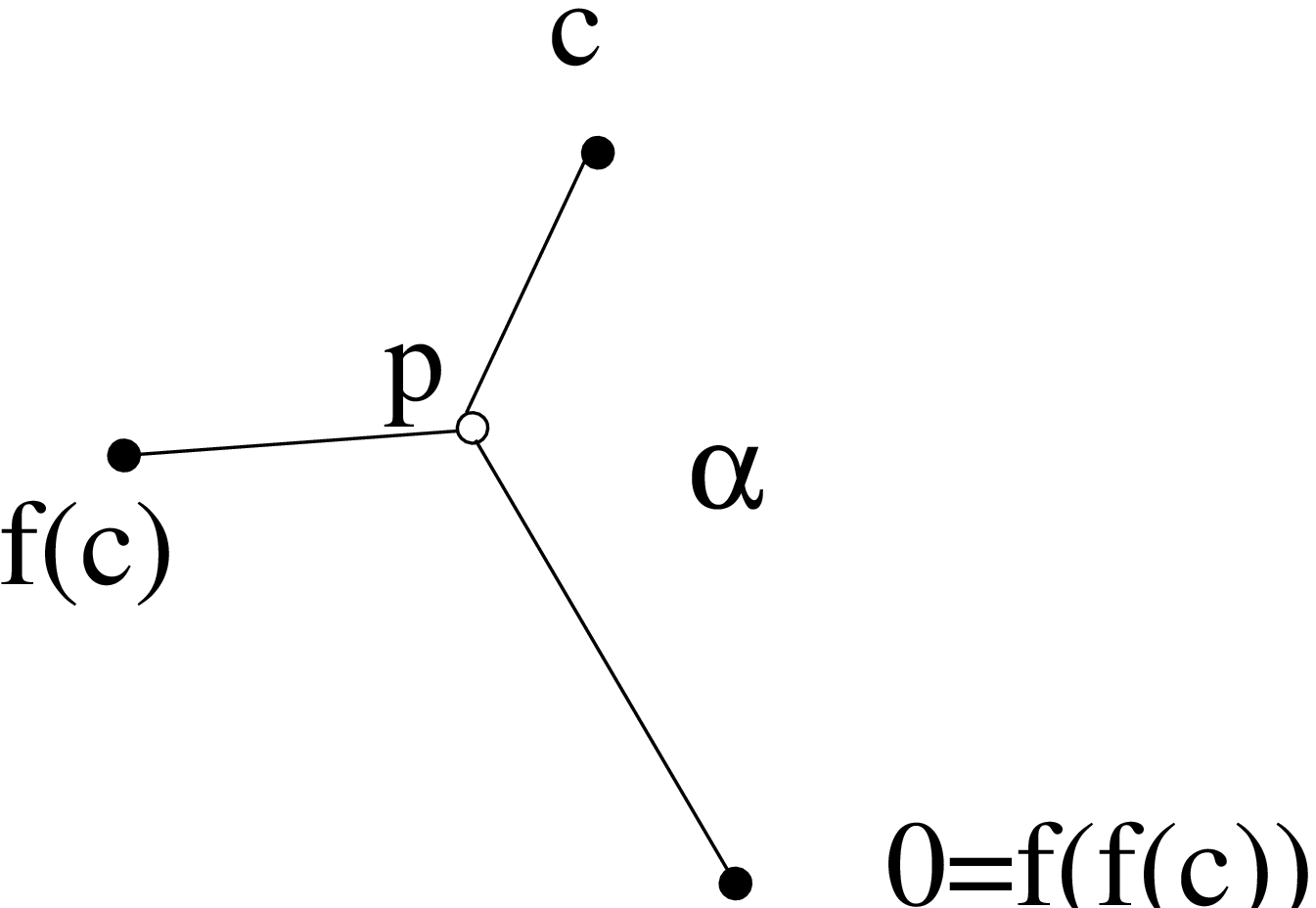}{1}{The tree for the Douady Rabbit.} 

We let $\alpha$ be the union of the two internal rays joining $0$ and
$c$ to $p$.  By Corollary 2, blowing up once along $\alpha$ yields a
branched covering equivalent to a rational map.  The Julia set of $f$
is shown in Figure \ref{fig:rabbit.ps}, and that of $f$ blown up along
$\alpha$ is shown in Figure
\ref{fig:rabbit_blow.ps}.  

\yincludeps{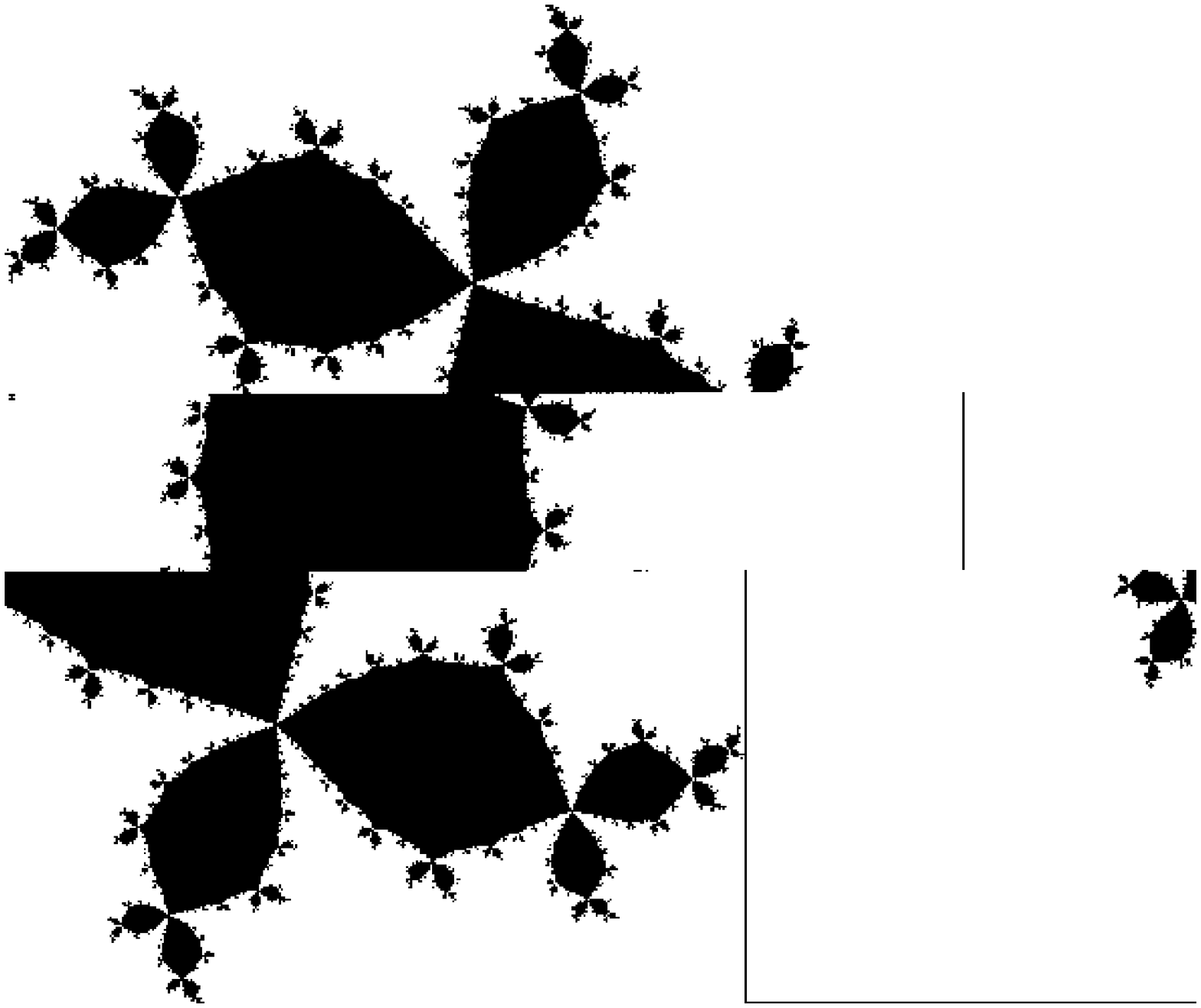}{1.75}{Douady's Rabbit.}
\yincludeps{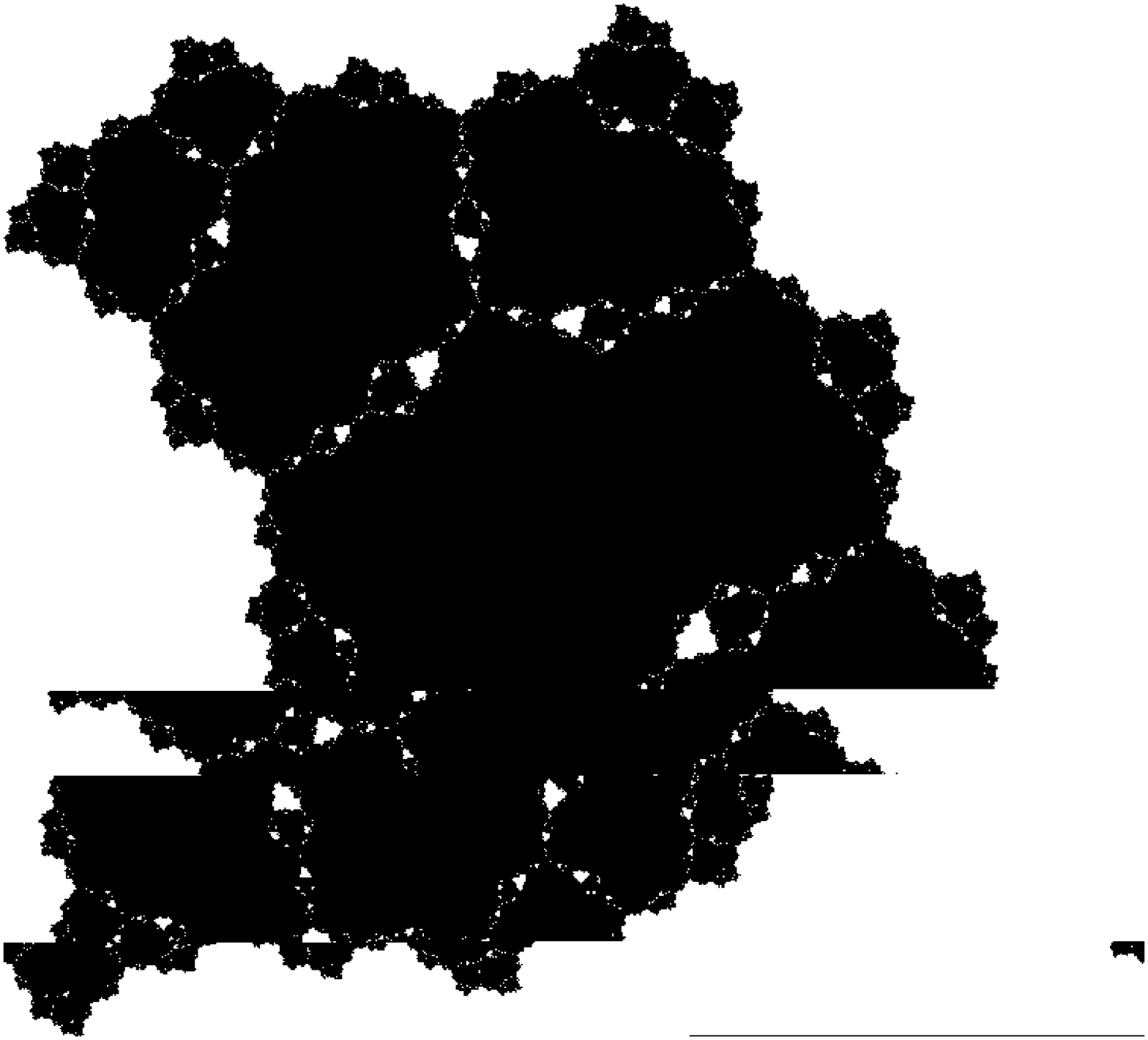}{1.75}{The Douady Rabbit blown up along $\alpha$.}

\smallbreak
\noindent{\bf Blowing up the ``airplane''.} Let $f(z)=z^{2}+c$, where
$c \approx -1.75488$ is chosen so that the critical point at zero is
periodic of period three with orbit lying in the real line.  Then $c <
0 < f(c)$.  Let $X=\{c,0,f(c)\}$.  The tree $T_{X}$ is 
the interval $[c,f(c)]$ with vertices $X$.  The polynomial $f$ sends
the invervals $[c,0]$ and $[0,f(c)]$ homeomorphically to the interval
$[c,f(c)]$. We let $\alpha=[c,0] \subset T_{X}$.  Then the 
Blowing Up conditions are satisfied, and by Theorem B, $(f,X)$ blown
up once along $\alpha$ is equivalent to a rational map.

The Julia set of $f$ is shown in Figure \ref{fig:airplane.ps}. 
The Julia set of the result of blowing up $f$ along $\alpha$ is shown
in Figure \ref{fig:airplane_blow.ps}; the map is given by 
$g(z)=a \frac{(z-1)^{3}}{1-3z}$, where $a \approx -2.37123$, normalized so
that $1 \mapsto 0$ by local degree three, $0
\mapsto -a$ by local degree two, and $-a \mapsto 1$ by local degree
one.

\xincludeps{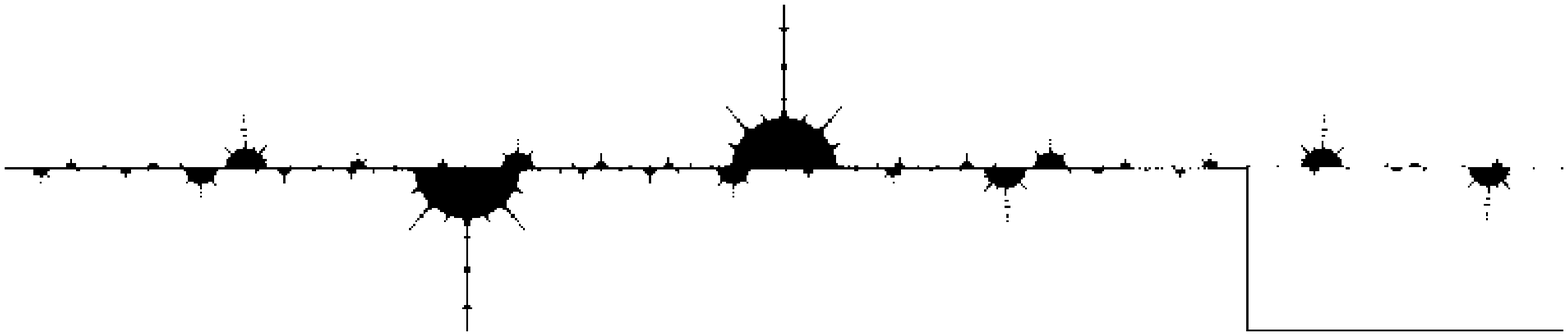}{4.5}{The airplane.}
\xincludeps{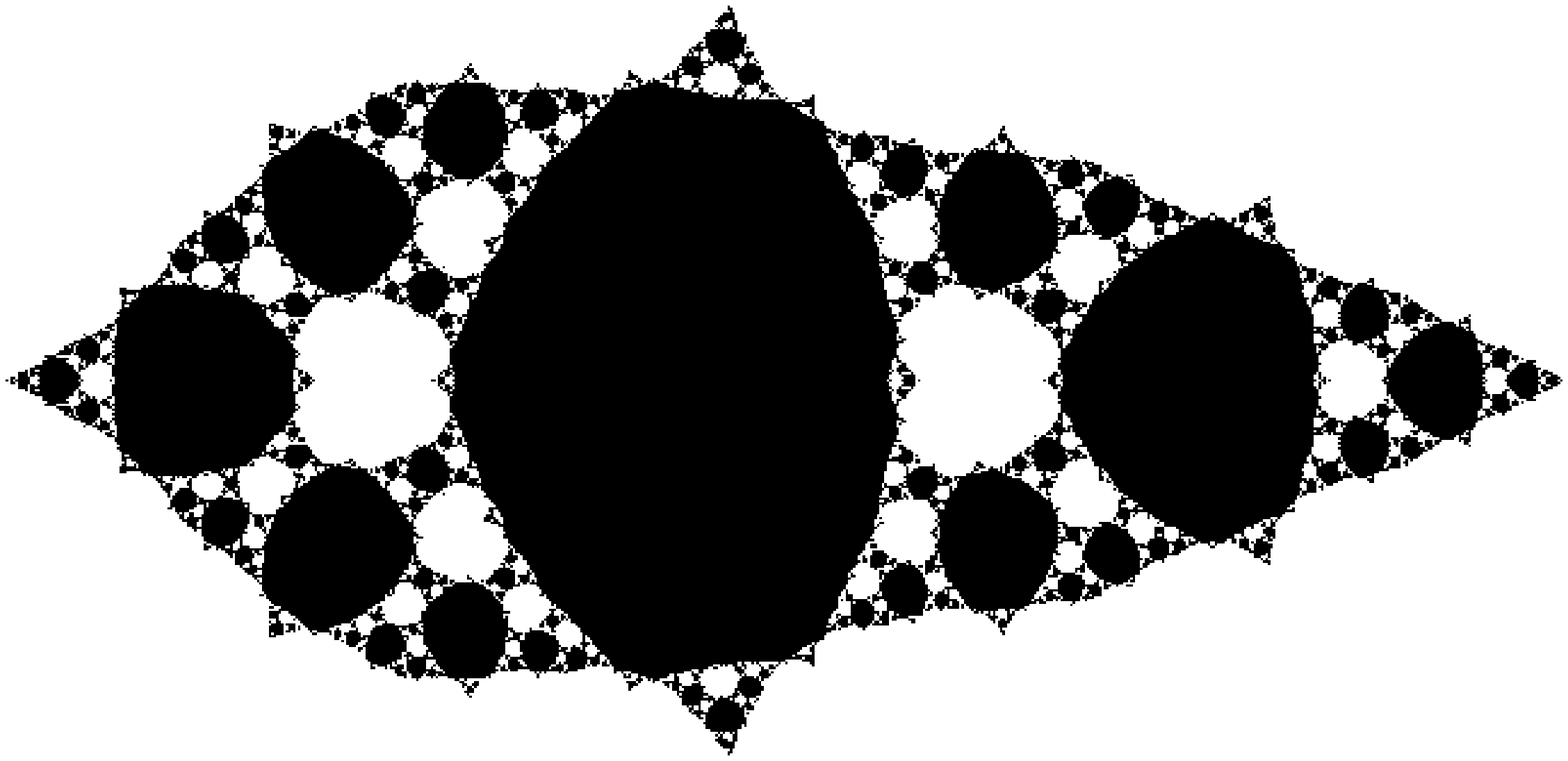}{4.5}{The airplane blown up along $\alpha$.}

\gap\gap

\noindent {\bf A degree four map with $S_3$ symmetry.} Let $G$ be a topological
graph in $S^2$ which is a triangle with vertices $X$ and edges $E$.
Let $M$ be the identity map.  By Theorem \ref{thm:blowing up moebius
transformations}, {\it Blowing up M\"obius transformations}, the
marked covering $(g,X)$ blown up once along each edge in $E$ is
equivalent to a degree four rational map $g$ with three fixed critical
points each of mulitiplicity two.  The full mapping class group of
orientation-preserving homeomorphisms $h$ of $S^{2}-X$ up to isotopy
is isomorphic to $S_3$ and commutes with $g$ up to equivalence since
the blowing up construction in this case is performed symmetrically.
That is, $h\circ g$ is equivalent to $g \circ h$.  It follows that the
conformal automorphism group of $g$ is $S_3$.  A formula for this map
is given in \cite{ctm:iterative} and is used there for giving a
generally convergent iterative algorithm for solving the cubic.

\gap\noindent{\bf Remark:} More generally one can apply Theorem 
\ref{thm:blowing up moebius transformations}, {\it
Blowing up M\"obius transformations} to other symmetric graphs and
the identity map to produce symmetric postcritically finite rational
maps; cf. \cite{ctm:pd:quintic} for a detailed discussion of rational
maps with icosahedral symmetry.

\section{Appendix: Proof of Theorem {\it Arcs intersecting
obstructions} }
\label{appendix:arcs}

%
%
\def\g{\gamma}
\def\G{\Gamma}
\def\la{\lambda}
\def\La{\Lambda}
\def\pf{postcritically finite }
\def\pfbc{postcritically finite branched covering }
\def\pfbcs{postcritically finite branched coverings }
\def\cc{connected component }
\def\ccs{connected components }
\def\gf{geometrically finite }
\def\pc{postcritical set }
\def\lc{locally connected }

\def\beginp{ {\noindent\em Proof.} }
\def\endp{\hspace*{\fill}$\rule{.55em}{.55em}$ \smallskip}

We will use the following easy fact:

\begin{lemma}
\label{lemma:crucial} Let $f$  be a  branched covering from $S^2$ to $S^2$. 
Let $B$ be a subset of $S^2$ such that $f|_B:B\to f(B)$ is a bijection
(resp. a $k$-to-$1$ mapping). Let $A$ be any subset of $S^2$. Then $\#
f^{-1}(A)\cap B=\# A\cap f(B)$ (resp. $\# f^{-1}(A)\cap B=k\cdot \#
A\cap f(B)$).
\end{lemma}

For $\La$ a curve system or an arc system in $(S^2,X))$, set
$\widetilde{\Lambda}$ (resp. $\widetilde{\Lambda}(f^n)$) to be the
union of those components of $f^{-1}(\Lambda)$ (resp. $f^{-n}(\La)$)
which are isotopic to elements of $\Lambda$. If $\La$ is irreducible,
each component of $\La$ is isotopic to some (possibly non-unique)
component of $\widetilde{\La}(f)$ (resp.  $\widetilde{\La}(f^n)$).

By the definition of irreducible arc system, there exists at least one
subset $\widetilde{\La}'$ of $\widetilde{\La}$ such that
$\widetilde{\La}'$ is isotopic to $\La$ and $f:\widetilde{\La}'-X \to
\La-X$ is a homeomorphism.  Then

\begin{equation}
\label{eqn:1} 
f_{\#,\G}(\g_j) \cdot \La
 \le \sum_{\g'\subset f^{-1}(\g_j)} \#(\g'\cap \widetilde{\La}')
=\#( f^{-1}(\g_j)\cap \widetilde{\La}') = \#(\g_j\cap \La) =
 \g_j\cdot \La\ 
\end{equation}
where the first inequality follows since $f:\widetilde{\La}'-X \to
\La-X$ is a homeomorphism and the second to last equality follows by
the lemma. 
Assume $f_{\#,\G}(\g_j)=\sum_i b_{ij}\g_i$. Then by (1), 
$ (\sum_j\sum_i b_{ij}\g_i) \cdot \La \le \G\cdot \La\ .$ 
But 
\begin{equation}
\label{eqn:2} (\sum_j\sum_i b_{ij}\g_i) \cdot \La =(\sum_i(\sum_j
b_{ij})\g_i) \cdot \La 
 \geq \sum_i \g_i \cdot \La = \G\cdot \La\ ,
\end{equation} 
where the last inequality is obtained as follows.  By irreducibility,
$f_\G$ has no zero row no zero column, so the same is true for
$f_{\#,\G}$.  But each entry of $f_{\#,\G}$ is a non negative integer,
so $b_{ij}\ge 1 $ if $b_{ij}\ne 0$.  We have then either A) $\G\cdot
\La=0$ or B) $\G\cdot
\La \neq 0$ and all inequalities above become equalities.

In case that $\G\cdot \La=0$, we claim that $\G\cdot f^{-1}(\La)=0 $.
Let $\g_i\in \G$. Then $\g_i\sim \g'$ with $\g'$ a curve in
$f^{-1}(\g_j)$ for some $\g_j\in \G$. Then $\g_i \cdot f^{-1}(\La) \le
\#(\g'\cap f^{-1}(\La))=k\#( \g_j\cap \La) = 0$.  The case for $n > 1$
is proved similarly.  This prove part one of the theorem. 

In case B), equalities in Equation \ref{eqn:2} imply that for each
$i$, there is a unique $j$ such that $b_{ij}\ne 0$ and that
$b_{ij}=1$. Thus for each $\g_i$, there is a unique $\g_j$ such that
some component of $f^{-1}(\g_j)$ is isotopic to $\g_i$. Also, there is
a unique such component $\g'$ of $f^{-1}(\g_j)$ isotopic to $\g_i$.
Since $f_{\#,\G}$ is irreducible, this implies that $f_{\#,\G}$ is a
transitive permutation of the basis vectors. As a consequence, $\g'$
is also the unique curve in $f^{-1}(\g_j)$ isotopic to an element of
$\G$.

Equalites in Equation \ref{eqn:1} tells us that $f_{\#,\G}(\g_j)\cdot \La
=\#(f^{-1}(\g_j)\cap \widetilde{\La}') =\g_j\cdot \La$.  Combining
with the facts that $f_{\#,\G}(\g_j) =[\g']_\G=\g_i$, $\g'\cdot \La\le
\#(\g'\cap \widetilde{\La}')$ and $\g'\cap \widetilde{\La}'\subset
f^{-1}(\g_j)\cap \widetilde{\La}'$, we get $\g_i\cdot \La=\g_j\cdot
\La$ and $ \g'\cap \widetilde{\La}' = f^{-1}(\g_j)\cap
\widetilde{\La}'$.  As a consequence, $\g'$ is also the unique curve
in $f^{-1}(\g_j)$ such that $\g'\cap \widetilde{\La}'\ne
\emptyset$. Moreover $f^{-1}(\G)\cap \widetilde{\La}'=\widetilde{\G}\cap
\widetilde{\La}'$.

The matrix $f_\G$ is a permutation with entry $a_{ij}=1/k_j$, where
 $k_j$ is the degree of $f:\g'\to \g_j$. Since $\la(f_\G)\geq 1$, we
 have $k_j=1$ for all $j$ and $f_\G=f_{\#,\G}$.  Thus
 $f:\widetilde{\G}\to \G$ is a homeomorphism.

Note that in the above argument we only used the following properties
of the triple $(\G,\La,\widetilde{\La}')$: the matrix $f_{\G}$ has
neither zero row nor zero column, with leading eigenvalue at least
$1$, the set $\widetilde{\La}'$ is isotopic to $\La$ and
$f:\widetilde{\La}'-X\to \La-X$ is a homeomorphism. By what we just
proved, the triple $(\La,\G,\widetilde{\G})$ verifies also these
properties. So we can redo the same argument for this new triple. We
conclude then for each $\la\in \La$, there is exactly one curve $\la'$
in $f^{-1}(\la)$ with $\la'\cap \widetilde{\G}\ne \emptyset$.
Moreover, $\la'$ is the unique curve in $f^{-1}(\la)$ which is
isotopic to an element of $\La$. The matrix $f_{\#,\La}$ is a
transitive permutation. As a consequence,
$\widetilde{\La}'=\widetilde{\La}$.

Replace $\widetilde{\La}'$ by $\widetilde{\La}$ in the beginning of
the proof, we conclude finally that $\widetilde{\G}\cap
f^{-1}(\La)=\widetilde{\G}\cap \widetilde{\La} = f^{-1}(\G)\cap
\widetilde{\La}$ and $\#\widetilde{\G}\cap \widetilde{\La}=\#\G\cap
\La=\G\cdot \La$. As a consequence, $\widetilde{\G}\cap (f^{-1}(\La)-
\widetilde{\La})=\emptyset$ and $(f^{-1}(\G)-\widetilde{\G})\cap
\widetilde{\La}=\emptyset$. 

Now take any arc $\la$ isotopic to some arc in $f^{-1}(\La)$ but to
none of the arcs in $\La$, it is isotopic to an arc $\la'\in
f^{-1}(\La)-\widetilde{\La}$. So $\G\cdot \la\le \#(\widetilde{\G}\cap
\la')=0$.

This proves the part $2(a),2(b)$ and $(2(d),n=1)$ of the theorem. To
get $2(c)$ and $(2(d),n>1)$, we just need to note that if $\G,\La$ are
irreducible for $f$, they are also irreducible for $f^n$, for all
$n\ge 1$. This gives the theorem.
\qed

The following theorem can be proved in exactly the same manner
as Theorem \ref{thm:arcs intersecting obstructions}.  

\begin{thm}[Arcs intersecting arcs]
Let $\Lambda_{1}, \Lambda_{2}$ be two irreducible systems of arcs.
Then either

1. $\Lambda_{1} \cdot \Lambda_{2} = 0$ and $ f^{-n}(\Lambda_{1}) \cdot
f^{-n}(\Lambda_{2}) = 0$ for all $n\ge 1$, or

2. $\Lambda_{1} \cdot \Lambda_{2} \neq 0$ and $(
f^{-n}(\Lambda_{1})-\widetilde{\Lambda}_{1}(f^n) ) \cdot \La_2=0=
\La_1\cdot ( f^{-n}(\Lambda_{2}) -\widetilde{ \Lambda}_{2}(f^n))$ for
all $n\ge 1$.  In this case we also have (a) The transformations
$f_{\sharp ,\Lambda_{i}}$ are transitive permutations of the basis
vectors; (b) for each $\lambda \in \Lambda_{i}$, there is a unique
component of $f^{-1}(\lambda)$ intersecting $\Lambda_{i'}$, $i \neq i'
\in \{1,2\}$; (c) (a) and (b) remain true if we replace $f$ by $f^n, n
\geq 1$. 
\end{thm}

\begin{cor}
If $\Lambda_{1}$ contains more than one periodic cycle, or if some
element of $\Lambda$ maps with multiplicity greater than one, then
$\Lambda_{1} \cdot \Lambda_{2} = 0$ for any irreducible arc system
$\Lambda_{2}$.  
\end{cor}


\end{document}